\documentclass[final,5p,times,twocolumn]{elsarticle}

\usepackage{amssymb}
\usepackage{amsmath}
\usepackage{bbm}
\usepackage{graphicx}
\usepackage{exscale}

\biboptions{sort&compress}

\usepackage{clrscode}
\usepackage{array}
\usepackage{hhline}
\usepackage{color, colortbl}
\definecolor{myGray}{rgb}{0.95,0.95,0.95}

\newcommand{\E}[1]{\mathop{{\rm \bf E}\!\left\{#1\right\}}\nolimits}

\newtheorem{lmm}{Lemma}
\newproof{pf}{Proof}

\newdefinition{rmk}{Remark}
\newdefinition{dfn}{Definition}
\newtheorem{exmp}{Example}

\begin{document}

\begin{frontmatter}
\title{On derivative-free extended Kalman filtering and its Matlab-oriented square-root implementations for state estimation in continuous-discrete nonlinear stochastic systems}

\author[CEMAT]{Maria V. Kulikova\corref{cor}} \ead{maria.kulikova@ist.utl.pt} \cortext[cor]{Corresponding
author.}

\author[CEMAT]{Gennady Yu. Kulikov} \ead{gennady.kulikov@tecnico.ulisboa.pt}

\address[CEMAT]{CEMAT, Instituto Superior T\'ecnico, Universidade de Lisboa, Av.~Rovisco Pais, 1049-001 Lisboa, Portugal.}

\begin{abstract}
Recent research in nonlinear filtering and signal processing has suggested an efficient derivative-free Extended Kalman filter (EKF) designed for discrete-time stochastic systems. Such approach, however, has failed to address the estimation problem for continuous-discrete models. In this paper, we develop a novel continuous-discrete derivative-free EKF methodology by deriving the related moment differential equations (MDEs) and sample point differential equations (SPDEs). Additionally, we derive their Cholesky-based square-root MDEs and SPDEs and obtain several numerically stable derivative-free EKF methods. Finally, we propose the MATLAB-oriented implementations for all continuous-discrete derivative-free EKF algorithms derived. They are easy to implement because of the built-in fashion of the MATLAB numerical integrators utilized for solving either the MDEs or SPDEs in use, which are the ordinary differential equations (ODEs). More importantly, these are accurate derivative-free EKF implementations because any built-in  MATLAB ODE solver includes the discretization error control that bounds the discretization error arisen and makes the implementation methods accurate. Besides, this is done in automatic way and no extra coding is required from users. The new filters are particularly effective for working with stochastic systems with highly nonlinear and/or nondifferentiable drift and observation functions, i.e. when the calculation of Jacobian matrices are either problematical or questionable. The performance of the novel filtering methods is demonstrated on several numerical tests.
\end{abstract}


\end{frontmatter}

\section{Introduction}\label{sect1}

The Extended Kalman filtering (EKF) framework is a widely used technique for solving the state estimation problem of nonlinear stochastic systems in practice. We may mention several applications that include air and space modeling with special attention to target tracking and orbit determination~\cite{CaSy83,LeMa82,MiHa99,PaBi97,PsHu93,TeSa08}, navigation and guidance of ships and underwater vehicles \cite{CaBr99,FoSa96,FuGr83}, road traffic positioning and estimation \cite{HoMu00,Ra97,ThWa99}, manipulator and robot control~\cite{LeDu91,MuKl97,NiNy99}, mobile station localization \cite{Sp99}, data assimilation in meteorology and related sciences \cite{Ev92,GhMa91,LeSt78,RoVa98,VoHe99}, industrial research \cite{AlAc97,Chen95,CrCh96,DiAu17,GuAu22,MuBu01,PrIr00,Wilson1998}, biology \cite{BoPe00,EnSa98,HiKa95}, speech recognition \cite{PaDo00}, image calibration \cite{CaNi98}, financial industry \cite{BhCh96,RiPe98}, etc. The EKF methodology is simple enough for implementation and it provides good estimation quality. These facts explain its popularity for solving practical problems. Despite the mentioned attractive features of the EKF, one of the main disadvantages is its requirement to calculate Jacobian matrices of the drift and measurement vector-functions of nonlinear model at hand. Recently, this problem has been resolved by the {\it derivative-free} EKF methodology suggested for estimating the state of {\it discrete-time} stochastic systems in~\cite{quine2006derivative}. The cited EKF-type framework is particularly effective for working with stochastic systems with highly nonlinear and/or nondifferentiable drift and observation functions, i.e. when the calculation of Jacobian matrices are either problematical or questionable. Such approach, however, has failed to address the estimation problem for a hidden state of continuous-discrete systems. The goal of our research is to fill in this gap, i.e. to derive the {\it continuous-discrete} derivative-free EKF framework.

It is well known in engineering literature that two basic approaches are routinely adopted to the development of {\it continuous-discrete} KF-like nonlinear filtering methods~\cite{Ja70,Frog2012,KuKu14IEEE_TAC}. The first one implies the use of numerical schemes for solving the given {\it stochastic differential equation} (SDE) of the system at hand. It is usually done by using either the Euler-Maruyama method or higher order methods, for example, that based on the It\^{o}-Taylor expansion~\cite{KlPl99}. An alternative implementation framework assumes the derivation of the related filters' moment differential equations (MDEs) and then an utilization of the numerical methods derived for solving ordinary differential equations (ODEs). We refer to the detailed discussion of advantages and disadvantages of these two derivation frameworks applied to continuous-discrete Unscented Kalman filter in~\cite{KuKu22Automatica}.  In the cited paper, it is shown that the approach based on the MDEs derivation yields accurate filters' implementation methods, in general. It is explained by the fact that the modern ODE solvers include the discretization error control, i.e. the error arisen is possible to control and to bound and this yields the accurate implementation way. Furthermore, the modern ODE solvers are self-adaptive algorithms, i.e. they generate the {\it adaptive} integration mesh in automatic way depending on the discretization error control involved for keeping the error within the prefixed tolerance. Such smart implementation strategy provides the accurate, self-adaptive and flexible filters' implementation methods where the users are requested to fix the tolerance value, only. They are useful for solving the state estimation problems with irregular sampling intervals (e.g., in case of missing measurements) as well as for working with stiff systems. This motivates us to derive the novel {\it continuous-discrete} derivative-free EKF methods within the discussed approach.

It is worth noting here that the {\it derivative-free} EKF principle introduced in~\cite{quine2006derivative} for discrete-time systems assumes a possibility to derive the sample point differential equations (SPDEs) as well. They can be used at the prediction filtering step instead of the MDEs to be solved. In our research we derive two implementation frameworks for continuous-discrete derivative-free EKF and discuss each of them in detail. To summarize, the contribution of this paper is threefold: (i) using the derivative-free EKF principle in~\cite{quine2006derivative}, we further derive the related MDEs and SPDEs for the derivative-free EKF prediction step and obtain the conventional MDE- and SPDE-based derivative-free EKF implementation techniques, respectively; (2) we next derive the square-root version of the MDEs for the Cholesky factors propagation as well as several square-root implementations for the measurement update step and obtain the numerically stable Cholesky-based square-root MDE- and SPDE-based derivative-free EKF methods; (3) we propose the MATLAB-oriented implementations for all continuous-discrete derivative-free EKF algorithms derived. They are easy to implement because of the built-in fashion of the MATLAB numerical ODEs integrators utilized for solving either the MDEs or SPDEs in use. More importantly, these are accurate derivative-free EKF implementations because any built-in  MATLAB ODE solver includes the discretization error control that bounds the discretization error arisen and makes the implementation methods accurate. Besides, this is done in automatic way and no extra coding is required from users. Finally, the performance of the novel filtering methods is demonstrated on several numerical tests.

The flow of the paper is as follows. Section~\ref{sect2} explains the derivative-free EKF principle and summarizes the related mean and covariance approximations at the measurement update step, which is the same for novel continuous- and previously published discrete-time derivative-free EKF estimator proposed in~\cite{quine2006derivative}. Section~\ref{sec:main:conventional} includes one of the main results of current research that is the derivation of the MDE and SPDE formulas. In this section, we also propose the conventional MATLAB-oriented implementation methods for the novel continuous-discrete derivative-free EKF derived. In Section~\ref{sec:main:squareroot}, the Cholesky-based square-root derivative-free EKF methods are developed. For that, the derivative-free MDEs are re-derived in terms of Cholesky factors of the filter covariance matrix. The results of numerical tests and comparative study with the standard continuous-discrete EKF are presented in Section~\ref{numerical:experiments}. Finally, Section~\ref{Section:conclusion} summarizes the key findings of our study as well as  some open questions for a future research.

\section{The derivative-free EKF-based principle: The mean and covariance approximation at the measurement update}\label{sect2}

Following~\cite{Haykin2009}, the Gaussianity assumption underlying the state estimation methods is a standard one utilized for the state estimation methods' derivation, e.g., the EKF, the cubature KF (CKF), the Gauss-Hermite KF (GHKF) filters and so on. This implies that all such filtering solutions are based on the first two moment calculations, i.e. the expectation and covariance of the prior and posterior random distributions should be computed. In this paper, we explore the derivative-free EKF framework introduced in~\cite{quine2006derivative} and briefly explain the related mean and covariance calculation formulas derived for  the measurement update step of the derivative-free EKF proposed in the cited paper.  The measurement processing step of the novel continuous-discrete derivative-free EKF technique coincides with the same step in the discrete-time derivative-free EKF estimator.

Following~\cite{quine2006derivative}, the key idea behind the derivative-free EKF estimation approach is an approximation of the probabilistic description of the state estimate by a minimally required number of a discrete vector set. More precisely, $n$ linearly independent vectors are demanded to span the $n$-dimensional vector space. Recall, at the measurement update step at time instance $t_k$, the one-step ahead predicted estimate together with the error covariance matrix are calculated on the previous filter's iterate, i.e. we have
\begin{align*}
\hat x_{k|k-1} & := \hat x(t_k|t_{k-1}) = \E{x(t_k)|z_1 \ldots z_{k-1}} \\
P_{k|k-1}      & := P(t_k|t_{k-1}) = \E{[x(t_k) - \hat x_{k|k-1}][ x(t_k) - \hat x_{k|k-1} ]^{\top}|z_1 \ldots z_{k-1}}
\end{align*}
where $\{z_1 \ldots z_{k-1}\}$ is the data history available for the filter from the measurement device, which is related to the hidden state vector through the measurement model given by
\begin{align}
z_k   & =  h(k,x(t_{k}))+v_k, \quad k =1,2,\ldots \label{eq1.2}
\end{align}
where  $x(t_k)$ is the $n$-dimensional unknown state vector to be estimated and $h:\mathbb R\times\mathbb
R^{n}\to\mathbb R^{m}$ is the observation function. The $m$-dimensional measurement vector $z_k = z(t_{k})$ comes at some discrete-time points $t_k$ with the sampling rate (sampling period) $\Delta_k=t_{k}-t_{k-1}$. The measurement noise term $v_k$ is assumed to be a white Gaussian noise with the zero mean and known covariance $R_k>0$, $R_k \in \mathbb R^{m\times m}$.

Following the derivative-free EKF principle  suggested in~\cite{quine2006derivative}, we consider a set of $n$ discrete vectors $\xi_i$ that approximate the second moment, i.e. the filter error covariance $P_{k|k-1}$, as follows:
\begin{align}
P_{k|k-1} & = \frac{1}{n} \sum \limits_{i=1}^{n} \xi_i \xi_i^{\top}. \label{eq:vectors}
\end{align}

Let us define the matrix $[\xi_1 | \ldots | \xi_n]$ collected from vectors $\xi_i$ in such a way that they are located by columns in the discussed matrix. Thus, we may express formula~\eqref{eq:vectors} as follows:
\begin{align}
n P_{k|k-1} & = \Bigl[\xi_1 \;|\; \ldots \;|\; \xi_n\Bigr] \; \Bigl[\xi_1 \;|\; \ldots \;|\; \xi_n\Bigr]^{\top}. \label{eq:vectors:1}
\end{align}

Taking into account the similarity between the calculations in~\eqref{eq:vectors:1} and the matrix square-root\footnote{Throughout the paper, we consider the lower triangular Cholesky factors, i.e. $P = SS^{\top}$ where $S$ is a lower triangular matrix with positive diagonal elements~\cite[Section~1.3.2]{Bj15}.} multiplication $P_{k|k-1} = P_{k|k-1}^{1/2}P_{k|k-1}^{\top/2}$, one may conclude
\begin{align}
  \Bigl[\xi_1 \;|\; \ldots \;|\; \xi_n\Bigr] & = \sqrt{n} P_{k|k-1}^{1/2}.  \label{eq:vectors:2}
  \end{align}

Using these $n$-vectors, $\xi_i$, $i=1,\ldots,n$, we define a further set of vectors about the predicted state with the distribution scaled by a constant $\alpha$ as follows:
\begin{align}
{\mathcal X}_{i, k|k-1} & =  \hat x_{k|k-1} + \frac{1}{\alpha} \xi_i,  & i & = 1, \ldots, n \label{Yvectors:new}
\end{align}
where the scalar parameter $\alpha$ is used to scale the discrete distribution about the mean.

The error covariance $P_{k|k-1}$ might be recovered from vectors ${\mathcal X}_{i, k|k-1}$, $i = 1, \ldots, n$, by taking into account formula~\eqref{Yvectors:new} and approximation in~\eqref{eq:vectors} in the following way:
\begin{align}
P_{k|k-1} & = \frac{1}{n} \sum \limits_{i=1}^{n} \xi_i \xi_i^{\top} = \frac{\alpha^2}{n} \sum \limits_{i=1}^{n} [{\mathcal X}_{i, k|k-1} - \hat x_{k|k-1}][{\mathcal X}_{i, k|k-1} - \hat x_{k|k-1}]^{\top} \nonumber \\
& = \alpha^2 \E{[{\mathcal X}_{i, k|k-1} - \hat x_{k|k-1}][{\mathcal X}_{i, k|k-1} - \hat x_{k|k-1}]^{\top}}.
\end{align}

Similar idea is used at the measurement update step for the residual covariance matrix approximation. More precisely, the vectors ${\mathcal X}_{i, k|k-1}$, $i=1, \ldots,n$ in~\eqref{Yvectors:new} and the predicted state vector $\hat x_{k|k-1}$ are then propagated through the observation function $h(\cdot)$ to generate
\begin{align}
Z_{i, k|k-1} & = h\left(k,{\mathcal X}_{i, k|k-1}\right), & i & = 1, \ldots, n, \label{eq:Zvectors} \\
\hat z_{k|k-1} & = h\left(k,\hat x_{k|k-1}\right). \label{eq:Zmean}
\end{align}

Thus, the residual covariance matrix is calculated by using vectors $\hat z_{k|k-1}$ and $Z_{i, k|k-1}$, $i = 1, \ldots, n$, as follows:
\begin{align}
R_{e,k} & =\E{[z_k - \hat z_{k|k-1}][z_k - \hat z_{k|k-1}]^{\top}} \nonumber \\
& \approx \alpha^2 \E{ \bigl[Z_{i, k|k-1} -  \hat z_{k|k-1} \bigr] \bigl[Z_{i, k|k-1} -  \hat z_{k|k-1}\bigr]^\top}\nonumber \\
& = \frac{\alpha^2}{n}\sum_{i=1}^{n} \bigl[Z_{i, k|k-1} -  \hat z_{k|k-1} \bigr] \bigl[Z_{i, k|k-1} -  \hat z_{k|k-1}\bigr]^\top + R_k, \label{eq:approx:Rek}
\end{align}

A cross-covariance $P_{xz,k}$ may also be evaluated by
\begin{align}
P_{xz,k} & =\E{[x(t_k) - \hat x_{k|k-1}][z_k - \hat z_{k|k-1}]^{\top}} \nonumber \\
& \approx \frac{\alpha^2}{n}\sum_{i=1}^{n} \bigl[{\mathcal X}_{i, k|k-1} -  \hat x_{k|k-1} \bigr] \bigl[Z_{i, k|k-1} -  \hat z_{k|k-1}\bigr]^\top. \label{eq:approx:Pxz}
\end{align}

Finally, the filtered estimate and the error covariance matrix of the derivative-free EKF are calculated as in the standard EKF, i.e. by the formulas (also see the algorithm in~\cite[Section~5]{quine2006derivative}):
\begin{align}
\hat x_{k|k} & =\hat x_{k|k-1}+{K}_k(z_k-\hat z_{k|k-1}),  \label{ckf:state}\\
{K}_{k} & = P_{xz,k}R_{e,k}^{-1}, \mbox{ and } P_{k|k}  = P_{k|k-1} - {K}_k R_{e,k} {K}_k^{\top}. \label{ckf:gain}
\end{align}

It is proved in~\cite{quine2006derivative} that in the limiting case as $\alpha \to \infty$, the derivative-free EKF formulas~\eqref{eq:approx:Rek}, \eqref{eq:approx:Pxz} converge to the standard EKF equations for calculating the  residual covariance $R_{e,k}$ and cross-covariance $P_{xz,k}$, respectively. The choice of the parameter $\alpha$ is adaptable and can be tuned on the application model at hand. More details can be found in the cited paper. The simulation results provided there show an excellent convergence to the standard EKF when $\alpha = 1000$.

The main goal of this paper is to derive the continuous-discrete derivative-free EKF methods and their useful MATLAB implementations. Following~\cite{higham2002nine}, any MATLAB implementation can be improved by vectorizing the operations, i.e. by maximizing the use of MATLAB's built-in functions that operate directly on vectors and matrices and by avoiding the loops. This feature dramatically reduces the execution times in MATLAB by using high-level operations on arrays rather than computing with individual components. In our case, the calculations by formulas~\eqref{eq:approx:Rek} and~\eqref{eq:approx:Pxz} should be vectorized.

First, let us denote the matrix ${\mathbb  X}_{k|k-1}$ collected from vectors ${\mathcal X}_{i, k|k-1}$, $i=1,\ldots, n$, defined in~\eqref{Yvectors:new} and located by columns, i.e.
\begin{align}
{\mathbb  X}_{k|k-1} & =\Bigl[{\mathcal X}_{1,k|k-1} \;|\; \ldots \;|\; {\mathcal X}_{n,k|k-1}\Bigr]. \label{notation:Xmatrix}
\end{align}

Taking into account formula~\eqref{eq:vectors:2}, it is not difficult to show that equation~\eqref{Yvectors:new} has the following simple matrix-vector form:
\begin{align}
{\mathbb  X}_{k|k-1}=\hat x_{k|k-1} {\mathbf 1}^{\top} + \frac{\sqrt{n}}{\alpha}P_{k|k-1}^{1/2} \label{newnewnewformula}
\end{align}
where ${\mathbf 1}$ is a column vector, whose elements are all equal to one.

Let us introduce the mean-adjusted and scaled matrix $\overline{\mathbb  X}_{k|k-1}$ obtained from the matrix ${\mathbb  X}_{k|k-1}$ calculated in~\eqref{notation:Xmatrix} as follows:
 \begin{align}
\overline{\mathbb  X}_{k|k-1} & = \frac{\alpha}{\sqrt{n}} \bigl[{\mathcal X}_{1,k|k-1} - \hat x_{k|k-1} | \ldots | {\mathcal X}_{n,k|k-1} - \hat x_{k|k-1}\bigr]. \label{notation:Xmatrix:center}
\end{align}
Similarly, we define the matrix $\overline{\mathbb  Z}_{k|k-1}$ collected from vectors ${Z}_{i, k|k-1}$, $i=1,\ldots, n$, in~\eqref{eq:Zvectors} in such a way that they are located by columns in the discussed matrix, i.e.
\begin{align}
\overline{\mathbb  Z}_{k|k-1} & = \frac{\alpha}{\sqrt{n}}\bigl[{Z}_{1,k|k-1} - \hat z_{k|k-1}| \ldots | {Z}_{n,k|k-1} - \hat z_{k|k-1}\bigr]. \label{notation:Zmatrix:center}
\end{align}

Taking into account~\eqref{notation:Xmatrix:center}, \eqref{notation:Zmatrix:center}, the derivative-free EKF formulas~\eqref{eq:approx:Rek}, \eqref{eq:approx:Pxz} can be calculated as follows:
\begin{align}
R_{e,k} & =  \overline{\mathbb  Z}_{k|k-1} \overline{\mathbb  Z}_{k|k-1}^\top + R_k, & P_{xz,k} & = \overline{\mathbb  X}_{k|k-1} \overline{\mathbb  Z}_{k|k-1}^\top. \label{eq:approx:Pxz:new}
\end{align}

In summary, the derivative EKF method has been derived for the discrete stochastic systems in~\cite{quine2006derivative}. The readers may refer to Section~5 in the cited paper to the implementation details of the {\it discrete-time} derivative-free EKF and its measurement update step discussed. In the next section, we explore a more sophisticated problem that is the derivation of the {\it continuous-time} derivative-free EKF technique.

\section{Derivation of continuous-discrete derivative-free EKF} \label{sec:main:conventional}

Consider continuous-discrete stochastic system with the measurement equation~\eqref{eq1.2} and the process equation given by
\begin{align}
dx(t) & = f\bigl(t,x(t)\bigr)dt+Gd\beta(t), \quad t>0,  \label{eq1.1}
\end{align}
where  $x(t)$ is the $n$-dimensional unknown state vector to be estimated and  $f:\mathbb R\times\mathbb
R^{n}\to\mathbb R^{n} $ is the time-variant drift function. The process uncertainty is represented by the additive noise term where
$G \in \mathbb R^{n\times q}$ is the time-invariant diffusion matrix and $\beta(t)$ is the $q$-dimensional Brownian motion whose increment $d\beta(t)$ is Gaussian white process independent of $x(t)$ and has the covariance $Q\,dt>0$. Recall, the $m$-dimensional measurement $z_k = z(t_{k})$ obeys equation~\eqref{eq1.2} and comes at some discrete-time points $t_k$ with the sampling rate (sampling period) $\Delta_k=t_{k}-t_{k-1}$. The measurement noise term $v_k$ in equation~\eqref{eq1.2} is assumed to be a white Gaussian noise with the zero mean and known covariance $R_k>0$, $R_k \in \mathbb R^{m\times m}$. Finally, the initial state $x(t_0)$ and the noise processes are assumed to be statistically independent, and $x(t_0) \sim {\mathcal N}(\bar x_0,\Pi_0)$, $\Pi_0 > 0$.

Let us explore the sampling interval $[t_{k-1}, t_{k}]$ and discuss the prediction step between two successive sampling time instants $t_{k-1}$ and $t_{k}$. Following~\cite[Example~4.19]{Ja70} and taking into account that no measurement has been received between $t_{k-1}$ and $t_{k}$, we come to the following {\em initial value problem} (IVP) for time-propagation of the expectation vector $\hat x(t)$ and covariance matrix $P(t)$ in the random state $x(t)$ that obeys equation~\eqref{eq1.1}:
\begin{eqnarray}\label{eq3.1a}
\hat x'(t) & = & \E{f\bigl(t,x(t)\bigr)},  \\ \nonumber
P'(t)      & = & \E{x(t)\,f^\top\!\bigl(t,x(t)\bigr)} -\hat x(t)\,\E{f^\top\!\bigl(t,x(t)\bigr)} \nonumber\\
           & + & \E{f\bigl(t,x(t)\bigr)\,x^\top\!(t)}-\E{f\bigl(t,x(t)\bigr)}\hat x^\top\!(t)\nonumber\\
           & + & G(t)\,Q(t)\,G^\top\!(t), \label{eq3.1b}\\
\hat x(t_{k-1}) & = & \hat x_{k-1|k-1},\quad P(t_{k-1})=P_{k-1|k-1}\label{eq3.1c}
\end{eqnarray}
where the prime denotes the differentiation in time, the initial vector $\hat x_{k-1|k-1}$ and matrix $P_{k-1|k-1}$ are the mean and covariance of the filtering solution computed at the sampling instant $t_{k-1}$.

It should stressed that the differential equations (\ref{eq3.1a}) and (\ref{eq3.1b}) are not conventional Ordinary Differential Equations (ODEs) since the operator $\E{\cdot}$ utilized on the right-hand sides of these equations is a multidimensional integral depending on all moments of the solution vector $x(t)$ to the SDE (\ref{eq1.1}) due to a nonlinearity of the drift function $f(\cdot)$. This also means that no one numerical method used for solving ODEs could be applied to treat the IVP in~\eqref{eq3.1a}~-- \eqref{eq3.1c}, e.g., see more explanation in~\cite[p.~137,~168]{Ja70}.

To overcome the mentioned obstacle, one needs to simplify the IVP (\ref{eq3.1a})--(\ref{eq3.1c}). For that we follow the derivative-free EKF framework explained in the previous section. In particular, the derivation of the mean equation is similar to the standard EKF technology. Assuming a sufficiently smooth drift function ${f}\bigl(t,x(t)\bigr)$ of SDE (\ref{eq1.1}), one expands it in the Taylor series about the estimate vector $\hat{x}(t)$ as follows:
\begin{equation}\label{eq2.2}
f\bigl(t, x(t)\bigr)=f\bigl(t, \hat x(t)\bigr)+\partial_{x}f\bigl(t, \hat x(t)\bigr)\bigl({x}(t)-\hat{x}(t)\bigr)+ HOT
\end{equation}
where the square matrix $\partial_{x}f\bigl(t, \hat x(t)\bigr)$ of size~$n$ denotes the partial derivative (that is, the conventional Jacobian) of the drift function $f (\cdot)$ and term $HOT$ stands to the remaining higher order terms. Truncating the expansion of $f (\cdot)$ to first order, we obtain
\begin{align}
 \E{f\bigl(t, x(t)\bigr)} &\approx \E{ f\bigl(t, \hat x(t)\bigr)+\partial_{x}f\bigl(t, \hat x(t)\bigr)\bigl({x}(t)-\hat{x}(t)\bigr) } \nonumber \\
& = \E{ f\bigl(t, \hat x(t)\bigr)} = f\bigl(t, \hat x(t)\bigr). \label{MDE:X}
\end{align}

Next, we need to derive the moment differential equation with respect to covariance $P(t)$. For that, the derivative-free EKF framework is applied, i.e. we chose a set of $n$ discrete vectors $\xi_i$ that approximate the second moment $P(t)$ that is
\begin{align}
P(t) = & \E{[x(t) - \hat x(t)][x(t) - \hat x(t)]^{\top}}  = \frac{1}{n} \sum \limits_{i=1}^{n} \xi_i \xi_i^{\top}. \label{eq:vectors:1}
\end{align}

Let us define the matrix $[\xi_1 \;|\; \ldots \;|\; \xi_n]$ collected from these vectors in such a way that they are located by columns in the discussed matrix. Again, taking into account the similarity between the calculations in~\eqref{eq:vectors} and the matrix square-root multiplication $P(t) = P^{1/2}(t)P^{\top/2}(t)$, we arrive to relationship $ \bigl[\xi_1 \;|\; \ldots \;|\; \xi_n\bigr] = \sqrt{n} P^{1/2}(t)$. We emphasize that $P^{1/2}(t)$ is a lower triangular matrix with positive diagonal elements and $P^{\top/2}(t)$ stands for its transpose.

Using these $n$-vectors, $\xi_i$, $i=1,\ldots,n$, we define a further set of vectors about the state estimate with the distribution scaled by the same constant $\alpha$ as follows:
\begin{align}
{\mathcal X}_{i}(t) & =  \hat x(t) + \frac{1}{\alpha} \xi_i,  & i & = 1, \ldots, n \label{Xvectors:new:1}
\end{align}
and, next, we propagate these vectors ${\mathcal X}_{i}(t)$, $i=1,\ldots,n$, through the drift function $f(\cdot)$ to generate the following matrix:
\begin{align}
{\mathbb F}{\mathbb X}(t) & = \left[f\left(t,{\mathcal X}_{1}(t)\right) \;|\; \ldots \;|\; f\left(t,{\mathcal X}_{n}(t)\right)\right].  \label{FX:matrix}
\end{align}

Again, we introduce the matrix ${\mathbb  X}(t)$  collected from these vectors ${\mathcal X}_{i}(t)$, $i=1,\ldots, n$, in such a way that they are located by columns in the discussed matrix, i.e. ${\mathbb  X} (t) =\bigl[{\mathcal X}_{1}(t)\;|\; \ldots \;|\; {\mathcal X}_{n}(t)\bigr]$ and taking into account formulas~\eqref{eq:vectors:1}, \eqref{Xvectors:new:1}, one gets
\begin{align}
{\mathbb  X}(t)=\hat x(t) \; {\mathbf 1}^{\top} + \frac{\sqrt{n}}{\alpha}P^{1/2}(t) \label{SP:equation:1}
\end{align}
and, hence, the centered matrix is
\begin{align}
\overline{{\mathbb X}}(t) & =  \left[{\mathcal X}_{1}(t) - \hat{x}(t) \;|\; \ldots \;|\; {\mathcal X}_{n}(t) - \hat{x}(t) \right] = \frac{\sqrt{n}}{\alpha}P^{1/2}(t).  \label{X:matrix:center}
\end{align}

Taking into account equation~\eqref{MDE:X} and~\eqref{X:matrix:center}, we approximate the terms on the right-hand side of equation~\eqref{eq3.1b} as follows:
\begin{align}
& \E{x(t)\,f^\top\!\bigl(t,x(t)\bigr)}-\hat x(t)\,\E{f^\top\!\bigl(t,x(t)\bigr)} \nonumber \\
& = \E{\bigl[x(t)-\hat{x}(t)\bigr] \bigl[f\bigl(t,x(t)\bigr)-\E{f\bigl(t,x(t)\bigr)}\bigr]^\top}\nonumber\\
&\approx \frac{\alpha^2}{n}\sum_{i=1}^{n}\bigl[{\mathcal X}_{i}(t) - \hat{x}(t)\bigr]\bigl[f\bigl(t, {\mathcal X}_{i}(t) \bigr)-f\bigl(t, \hat x(t)\bigr)\bigr]^\top \nonumber \\
& = \frac{\alpha^2}{n} \overline{{\mathbb X}}(t) \; \overline{{\mathbb F}{\mathbb X}}^{\top}(t) =  \frac{\alpha}{\sqrt{n}} P^{1/2}(t) \; \overline{{\mathbb F}{\mathbb X}}^{\top}(t) \label{eq:pr:1}
\end{align}
where the centered matrices are defined by
{\small
\begin{align}
\overline{{\mathbb F}{\mathbb X}}(t) & = \bigl[f\left(t,{\mathcal X}_{1}(t)\right) - f\bigl(t, \hat x(t)\bigr)\;| \ldots |\; f\left(t,{\mathcal X}_{n}(t)\right) - f\bigl(t, \hat x(t)\bigr) \bigr].  \label{FX:matrix:center}
\end{align}
}

Similarly, we calculate the following terms:
\begin{align}
& \E{f\bigl(t,x(t)\bigr)\,x^\top\!(t)} - \E{f\bigl(t,x(t)\bigr)}\hat x^\top\!(t) \nonumber \\
& =\E{\bigl[f\bigl(t,x(t)\bigr)-\E{f\bigl(t,x(t)\bigr)}\bigr]\bigl[x(t)-\hat{x}(t)\bigr]^\top}\nonumber \\
&\approx \frac{\alpha^2}{n}\sum_{i=1}^{n}\bigl[f\bigl(t, {\mathcal X}_{i}(t) \bigr)-f\bigl(t, \hat x(t)\bigr)\bigr]\bigl[{\mathcal X}_{i}(t) - \hat{x}(t)\bigr]^\top \nonumber \\
& = \frac{\alpha^2}{n} \overline{{\mathbb F}{\mathbb X}}(t) \; \overline{{\mathbb X}}^{\top}(t) =  \frac{\alpha}{\sqrt{n}} \overline{{\mathbb F}{\mathbb X}}(t) P^{\top/2}(t)  \label{eq:pr:2}
\end{align}

Having substituted formulas~\eqref{MDE:X} and~\eqref{eq:pr:1}, \eqref{eq:pr:2} into the IVP (\ref{eq3.1a})--(\ref{eq3.1c}), we obtain the {\it Moment Differential Equations} to be solved at the prediction step of the continuous-discrete derivative-free EKF:
\begin{eqnarray}
\hat{x}'(t) & = & {f}\bigl(t, \hat x(t)\bigr), \label{eq2.4a} \\
{P}'(t)     & = &  \frac{\alpha}{\sqrt{n}} P^{1/2}(t) \; \overline{{\mathbb F}{\mathbb X}}^{\top}(t)  + \frac{\alpha}{\sqrt{n}} \overline{{\mathbb F}{\mathbb X}}(t) P^{\top/2}(t) \nonumber \\
& + & G(t)\,Q(t)\,G^\top\!(t), \label{eq2.4b} \\
\hat x(t_{k-1}) & = & \hat x_{k-1|k-1},\quad P(t_{k-1})=P_{k-1|k-1}. \label{eq2.4c}
\end{eqnarray}

The main disadvantage of the MDE-based continuous-discrete derivative-free EKF derived above is the assumption that the covariance matrix  ${P}(t)$ coming out of numerical integration of the MDE in~\eqref{eq2.4b} is a positive definite matrix and, hence, its Cholesky factorization  always  exists over each sampling interval $t \in [t_{k-1}, t_k]$. Unfortunately, this condition is often violated in practice because of the numerical integration and roundoff errors involved. We may improve the continuous-discrete derivative-free EKF methodology by designing a more stable version where the {\it sample points} ${\mathcal X}_{i}(t)$, $i=1, \ldots, n$, are propagated at the prediction step instead of $P(t)$. This yields the {\it Sample Point Differential Equations} (SPDEs) to be solved instead of the MDEs at the prediction step of the continuous-discrete derivative-free EKF. To derive the SPDE-based continuous-discrete derivative-free EKF technique, we first need to prove the following result.

\begin{lmm}\label{Lm3-1}
Let the MDE in~\eqref{eq2.4b} possess the unique solution ${P}(t)$ in a sampling period $[t_{k-1},t_{k}]$ and has factorization $P(t) = P^{1/2}(t)P^{\top/2}(t)$  where $P^{1/2}(t)$ stands for the lower triangular Cholesky factor and $P^{\top/2}(t)$ denotes its transpose. If the square-root factor $P^{1/2}(t)$ is sufficiently smooth in the sampling interval $[t_{k-1},t_{k}]$ then it satisfies the Square-Root MDEs (SR-MDEs):
\begin{align}
 \bigl[P^{1/2}(t)\bigr]' & = P^{1/2}(t) \; \Phi\Bigl( P^{-1/2}(t) \; M(t) \; P^{-\top/2}(t)\Bigr) \label{eq3.8b}
\end{align}
where the matrix $M(t)$ of size~$n\times n$ stands for
\begin{align}
M(t) &:= \frac{\alpha}{\sqrt{n}} \left[ P^{1/2}(t) \; \overline{{\mathbb F}{\mathbb X}}^{\top}(t) + \overline{{\mathbb F}{\mathbb X}}(t) P^{\top/2}(t)\right] \nonumber \\
& + G(t)\,Q(t)\,G^\top\!(t) \label{eq3.9}
\end{align}
where the factor ${P}^{-\top/2}(t)$ refers to the transposed version of the inverse square-root factor ${P}^{-1/2}(t)$ and the operator $\Phi\bigl(A(t)\bigr)$ is a time-variant and matrix-valued one. It returns the lower triangular matrix defined by
\begin{equation} \label{eq3.10}
\Phi\bigl(A(t)\bigr):=\overline{L}(t) + 0.5\,{D}(t)
\end{equation}
where $A(t)$ is fragmented into strictly lower triangular $\overline{L}(t)$, diagonal $D(t)$ and strictly upper triangular $\overline{U}(t)$ parts, i.e. ${A}(t)=\overline{L}(t)+{D}(t)+\overline{U}(t)$.
\end{lmm}
\begin{pf} Taking into account factorization $P(t) = P^{1/2}(t)P^{\top/2}(t)$ and substituting it to the left-hand side of the MDE equation~\eqref{eq2.4b}, we obtain
\begin{align}
    & \bigl[P^{1/2}(t)\bigr]' P^{\top/2}(t) + P^{1/2}(t) \bigl[P^{\top/2}(t)\bigr]'  =  \frac{\alpha}{\sqrt{n}} P^{1/2}(t) \; \overline{{\mathbb F}{\mathbb X}}^{\top}(t)    \nonumber \\
& \qquad + \frac{\alpha}{\sqrt{n}} \overline{{\mathbb F}{\mathbb X}}(t) P^{\top/2}(t) + G(t)\,Q(t)\,G^\top\!(t), \label{prpr1}
\end{align}

Next, we left-multiply both sides of formula~\eqref{prpr1} by the inverse square-root factor ${P}^{-1/2}(t)$ and, then, right-multiply by its transpose ${P}^{-\top/2}(t)$. We obtain
\begin{align}
& P^{-1/2}(t) \bigl[P^{1/2}(t)\bigr]'  +  \bigl[P^{\top/2}(t)\bigr]'  P^{-\top/2}(t)  \nonumber \\
& \qquad =  P^{-1/2}(t) \; M(t) \; P^{-\top/2}(t) \label{prpr2}
\end{align}
where matrix $M(t)$ is defined in~\eqref{eq3.9}.

Let us consider the matrix equation~\eqref{prpr2}. The left-hand side in formula~\eqref{prpr2} is the sum of the lower triangular matrix and the upper triangular one. In addition, they have the same diagonal entries because of the property $\bigl[{P}^{\top/2}\!(t)\bigr]'=\bigl[{P}'(t)\bigr]^{\top/2}$. Besides, its right-hand side is a symmetric matrix. Thus, the lower triangular matrix standing on the left-hand side of the matrix equality~\eqref{prpr2} should be equal to the lower triangular part of its right-hand matrix yielded by means of formula (\ref{eq3.9}) and according to operator (\ref{eq3.10}). This completes the proof of Lemma~\ref{Lm3-1}.
 \qed
\end{pf}

Lemma~\ref{Lm3-1} can be used further to derive the SPDEs for the time update step of the continuous-discrete derivative-free EKF technique. More precisely, we prove the following result.

\begin{lmm}\label{Lm3-2}
Let the MDE in~\eqref{eq2.4b} possess the unique solution ${P}(t)$ in a sampling period $[t_{k-1},t_{k}]$ and has factorization $P(t) = P^{1/2}(t)P^{\top/2}(t)$  where $P^{1/2}(t)$ stands for the lower triangular Cholesky factor and $P^{\top/2}(t)$ denotes its transpose. If the square-root factor $P^{1/2}(t)$ is sufficiently smooth in the sampling interval $[t_{k-1},t_{k}]$ then the sample points in~\eqref{SP:equation:1} will be also smooth and satisfy the following  SPDEs expressed in a matrix form:
\begin{align}
{\mathbb  X}'(t) & = {f}\bigl(t, \hat x(t)\bigr) \; {\mathbf 1}^{\top} \nonumber \\
 & + \frac{\sqrt{n}}{\alpha} P^{1/2}(t) \; \Phi\Bigl( P^{-1/2}(t) \; M(t) \; P^{-\top/2}(t)\Bigr) \label{SPDEs:eqs}
\end{align}
where the matrix ${\mathbb  X}'(t)$ is collected from the sample vectors ${\mathcal X}'_{i}(t)$, $i=1,\ldots, n$, in such a way that they are located by columns in the discussed matrix, i.e. ${\mathbb  X}'(t) =\bigl[{\mathcal X}'_{1}(t)| \ldots | {\mathcal X}'_{n}(t)\bigr]$. The matrix $M(t)$ of size~$n\times n$ stands for the matrix defined in equation~\eqref{eq3.9}. The operator $\Phi\bigl(A(t)\bigr)$ is defined in formula~\eqref{eq3.10}.
\end{lmm}
\begin{pf}
Having differentiated the matrix equation~\eqref{SP:equation:1}, we get
\begin{align}
{\mathbb  X}'(t)=\hat x'(t) \; {\mathbf 1}^{\top} + \frac{\sqrt{n}}{\alpha}\left[P^{1/2}(t)\right]'. \label{prpr10}
\end{align}
Having substituted the formula for calculating the derivative of the square-root factor $\left[P^{1/2}(t)\right]'$ proved in Lemma~\ref{Lm3-1} and the right-hand side of the MDE formula in~\eqref{eq2.4a}, we arrive at equality~\eqref{SPDEs:eqs}. This completes the proof.
\qed \end{pf}

Based on Lemma~\ref{Lm3-1} and~\ref{Lm3-2}, we formulate the MDE- and SPDE-based continuous-discrete derivative-free EKF techniques derived. We summarize them in the form of pseudo-codes presented in Table~\ref{Tab:1}. It should be stressed that methods proposed in this paper are the general Matlab-oriented implementations. This means that any MATLAB's built-in ODE solver can be utilized in the implementation algorithms developed. We denote the function to be utilized by \verb"odesolver". The users are free to choose any method from~\cite[Table~12.1]{Higham2005} and substitute instead of \verb"odesolver".

As mentioned in Introduction, the main advantage of using MATLAB's built-in ODE solvers is that they allow to control the numerical integration accuracy in automatic way while solving the MDEs and SPDEs at any sampling interval $[t_{k-1}, t_k]$. Indeed, they include the local discretization error control that bounds the discretization errors occurred. More precisely, the solver creates the adaptive variable stepsize mesh in such a way that the discretization error arisen is less than the tolerance value pre-defined by user. This is done in automatic way by MATLAB's built-in functions and no extra coding is required from users except for setting the ODE solvers' options prior to filtering as follows:
\begin{equation}\label{eq2.136}\scriptsize
\texttt{options = odeset('AbsTol',LET,'RelTol',LET,'MaxStep',0.1)}
\end{equation}
where the parameters \texttt{AbsTol} and \texttt{RelTol} determine portions of the {\em absolute} and {\em relative} local error utilized in the built-in control mechanization, respectively, and $0.1$ limits the maximum step size $\tau^{\rm max}$ for numerical stability reasons. Formula~\eqref{eq2.136} implies that \texttt{AbsTol=RelTol=LET} where the parameter \texttt{LET} sets the requested local accuracy of numerical integration with the MATLAB code.

\begin{table*}[ht!]
{\small
\renewcommand{\arraystretch}{1.3}
\caption{The {\it conventional} continuous-discrete MATLAB-oriented derivative-free EKF implementation methods.} \label{Tab:1}
\centering
\begin{tabular}{l||l|l}
\hline
& \cellcolor{myGray} {\bf MDE-based derivative-free EKF: Algorithm~1} & \cellcolor{myGray} {\bf SPDE-based derivative-free EKF: Algorithm~2} \\
\hline
\hline
\textsc{Initialization:}  &  \multicolumn{2}{l}{0. Set the initial values $\hat x_{0|0} = \bar x_0$, $P_{0|0} = \Pi_0$ and derivative-free EKF parameter $\alpha$. Set the ODE solver's options in~\eqref{eq2.136}.  }\\
\hline
\textsc{Time} & 1. Form matrix $XP_{k-1|k-1} = [\hat x_{k-1|k-1}, P_{k-1|k-1}]$. & 1. $P_{k-1|k-1}=P_{k-1|k-1}^{1/2}P_{k-1|k-1}^{\top/2}$, get ${\mathbb  X}_{k-1|k-1}=\hat x_{k-1|k-1} {\mathbf 1}^{\top} + \frac{\sqrt{n}}{\alpha}P_{k-1|k-1}^{1/2}$. \\
\textsc{Update (TU):} & 2. Reshape $x^{(0)}_{k-1} = XP_{k-1|k-1}\verb"(:)"$. & 2. Collect matrix $A = [\hat x_{k-1|k-1}, {\mathbb  X}_{k-1|k-1}]$. Reshape $\widetilde{{\mathbb  X}}^{(0)}_{k-1|k-1} = A\verb"(:)"$. \\
& 3. Integrate $x_{k|k-1}\leftarrow \texttt{odesolver[MDEs},x^{(0)}_{k-1},[t_{k-1},t_k]]$. & 3. $\widetilde{{\mathbb  X}}_{k|k-1}\leftarrow \texttt{odesolver[SPDEs},\widetilde{{\mathbb  X}}^{(0)}_{k-1|k-1},[t_{k-1},t_k]]$. \\
& 4. $XP_{k|k-1} \leftarrow \texttt{reshape}(x_{k|k-1}^{\texttt{end}},n,n+1)$.  & 4. ${X}_{k|k-1} \leftarrow \texttt{reshape}(\widetilde{{\mathbb  X}}_{k|k-1}^{\texttt{end}},n,n+1)$. \\
& 5. Recover $\hat x_{k|k-1} = [XP_{k|k-1}]_1$. & 5. Recover $\hat x_{k|k-1} =[{X}_{k|k-1}]_1$ and sample points ${\mathbb X}_{k|k-1} = [X_{k|k-1}]_{2~:~n+1}$. \\
& 6. Recover $P_{k|k-1} = [XP_{k|k-1}]_{2~:~n+1}$. & 6. Recover $P^{1/2}_{k|k-1} = \frac{\alpha}{\sqrt{n}}\verb"tril"\bigl({\mathbb  X}_{k|k-1}-\hat x_{k|k-1}\; {\mathbf 1}^{\top}\bigr)$. \\
\hline
\textsc{Measurement} & 7a. Cholesky decomposition $P_{k|k-1}=P_{k|k-1}^{1/2}P_{k|k-1}^{\top/2}$ & $-$ Square-root $P^{1/2}_{k|k-1}$ is already available from TU. \\
\textsc{Update (MU):} & 7b. Define all sample points ${\mathbb  X}_{k|k-1}=\hat x_{k|k-1} {\mathbf 1}^{\top} + \frac{\sqrt{n}}{\alpha}P_{k|k-1}^{1/2}$. & $-$ The matrix ${\mathbb  X}_{k|k-1}$ is already available from TU. \\
\cline{2-3}
& \multicolumn{2}{l}{8. Propagate all sample points ${\mathcal X}_{i,k|k-1}$ from ${\mathbb  X}_{k|k-1}$: ${Z}_{i,k|k-1} = h\left(k,{\mathcal X}_{i,k|k-1}\right)$, $i=1,\ldots,n$, and the estimate $\hat z_{k|k-1} = h\left(k,\hat x_{k|k-1}\right)$. } \\
& \multicolumn{2}{l}{9. Define $\overline{\mathbb  X}_{k|k-1}  = \frac{\alpha}{\sqrt{n}} \bigl[{\mathcal X}_{1,k|k-1} - \hat x_{k|k-1} | \ldots | {\mathcal X}_{n,k|k-1} - \hat x_{k|k-1}\bigr]$ and $\overline{\mathbb  Z}_{k|k-1}  = \frac{\alpha}{\sqrt{n}}\bigl[{Z}_{1,k|k-1} - \hat z_{k|k-1}| \ldots | {Z}_{n,k|k-1} - \hat z_{k|k-1}\bigr]$.  } \\
& \multicolumn{2}{l}{10. Compute $R_{e,k}=\overline{\mathbb Z}_{k|k-1}\overline{\mathbb Z}_{k|k-1}^{\top}+R_k$, $P_{xz,k}=\overline{\mathbb X}_{k|k-1}\overline{\mathbb Z}_{k|k-1}^{\top}$, ${K}_{k}=P_{xz,k}R_{e,k}^{-1}$, $\hat x_{k|k}=\hat x_{k|k-1}+{K}_k(z_k-\hat z_{k|k-1})$, $P_{k|k}=P_{k|k-1} - {K}_k R_{e,k} {K}_k^{\top}$.} \\
\hline
\hline
Auxiliary & $[\tilde x(t)] \leftarrow \proc{MDEs}(x(t),t,n,\alpha,G,Q)$
&
$[\widetilde{{\mathbb  X}}(t)] \leftarrow \proc{SPDEs}(\widetilde{{\mathbb  X}}(t),t,n,\alpha,G,Q)$\\
Functions & Get matrix $X = \verb"reshape"(x(t),n,n+1)$; & Get matrix $X \leftarrow \verb"reshape"(\widetilde{{\mathbb  X}}(t),n,n+1)$;\\
& Recover  $\hat x(t) = [X]_1$ and $P(t) = [X]_{2~:~n+1}$;            & Recover $\hat x(t) = [X]_{1}$ and all sample points ${\mathbb X}(t) = [X]_{2~:~n+1}$;\\
& Cholesky decomposition: $P(t) = P^{1/2}(t)P^{\top/2}(t)$;  & Recover $P^{1/2}(t) = \frac{\alpha}{\sqrt{n}}\verb"tril"\bigl({\mathbb  X}(t)-\hat x(t)\; {\mathbf 1}^{\top}\bigr)$;\\
& Define all sample points by~\eqref{SP:equation:1};
& Propagate $d{\hat x}(t)/dt:= f\bigl(t,{\hat x}(t)\bigr)$ and all sample points to find $\overline{\mathbb FX}(t)$;\\
& Propagate all sample points to compute $\overline{\mathbb FX}(t)$; &  Find $M(t)$ by~\eqref{eq3.9} and split $M=\bar L + D + \bar U$. Find $\Phi(M) = \bar L + 0.5 D$; \\
& Compute the right-hand side of the MDEs in~\eqref{eq2.4a}, \eqref{eq2.4b}; & Compute the right-hand side of~\eqref{SPDEs:eqs}, i.e. get $d{\mathbb  X}(t)/dt$;\\
& Get extended matrix $\tilde X = [d\hat x(t)/dt, dP(t)/dt]$; & Collect extended matrix $A= [d\hat x(t)/dt, \; d{\mathbb  X}(t)/dt]$;\\
& Reshape into a vector form $\tilde x(t)=\tilde X\verb"(:)"$. & Reshape into a vector form $\widetilde{{\mathbb  X}}(t)=A\verb"(:)"$.\\
\hline
\end{tabular}
}
\end{table*}

It should be taken into account that the MATLAB's built-in ODE solvers are vector-functions and, hence, one should re-arrange both the MDEs in~\eqref{eq2.4a}, \eqref{eq2.4b} and SPDEs in~\eqref{SPDEs:eqs} in the form of unique vector of functions, which is to be sent to the ODE solver. The MATLAB built-in function \texttt{reshape} performs this operation. More precisely, \texttt{A(:)} returns a single column vector of size $M\times N$ collected from columns of the given array $A \in {\mathbb R}^{N\times M}$. It is the same as \texttt{reshape(A,M*N,1)}. In general, \texttt{reshape(A,M,N)} returns the $M$-by-$N$ matrix whose elements are taken columnwise from $A$. Next, the built-in function \texttt{tril(A)} extracts a lower triangular part from an array $A$, i.e. the result is a lower triangular matrix. Finally, the notation $\texttt{[A]}_i$ stands for the i-th column of any matrix $A$, meanwhile $\texttt{[A]}_{i:j}$ means the matrix collected from the columns of $A$ taken from the $i$th column up to the $j$th one.

Let us discuss the MDE- and SPDE-based Matlab-oriented continuous-discrete derivative-free EKF techniques summarized in Table~\ref{Tab:1}.
We remark that both Algorithm~1 and~2 are of {\it conventional}-type implementations because the filter error covariance matrix is propagated and updated. They are also of {\it covariance}-type methods, i.e. the error covariance matrix is processed instead of its inverse (which is called the information matrix) as it is done in {\it information}-type implementations. Although the methods treat the full covariance matrix, their Cholesky decomposition is required in each iterate to generate the sample points. This is one of the sources of the derivative-free EKFs instability. Indeed, the covariance matrix  ${P}(t)$ coming out of numerical integration of the MDE in~\eqref{eq2.4b} as well as after the measurement update step should be a positive definite matrix to ensure the existence of Cholesky factorization required. The implementation methods fail when this condition is violated due to roundoff and/or numerical integration errors. Likewise to CKF, UKF, GHQF and other sample-data KF-like estimators, the derivation of the square-root implementations is of special interest for the novel derivative-free EKF to avoid the source of the methods' failure at any filtering step due to infeasible matrix factorization. To avoid this problem, we will derive the square-root MDE- and SPDE-based derivative-free EKF implementations in the next section.

Following Table~\ref{Tab:1}, the SPDE-based derivative-free EKF implementation in Algorithm~2 does not require the Cholesky factorization at the prediction step compared to the MDE-based Algorithm~1. It comes from the fact that the sample vectors are propagated in Algorithm~2 and, hence, they are directly available at the next measurement update step. In other words, Algorithm~2 skips the matrix factorization at step~7. It makes Algorithm~2 more numerically stable than Algorithm~1. It is also worth noting here that the square-root matrix $P^{1/2}_{k|k-1}$ might be recovered from the propagated sample points according to formula~\eqref{newnewnewformula} and as shown on step~6 of Algorithm~2. Note, the factor $P_{k|k-1}^{1/2}$ is required at step~10 of Algorithm~2 for calculating the filter covariance $P_{k|k}$ through the predicted value $P_{k|k-1} = P_{k|k-1}^{1/2}P_{k|k-1}^{\top/2}$.

The auxiliary functions summarized in Table~\ref{Tab:1} are intended for computing the right-hand side functions in MDEs~\eqref{eq2.4a}, \eqref{eq2.4b} and SPDEs~\eqref{SPDEs:eqs}, respectively. They are presented in a vector form required by any built-in MATLAB ODEs integration scheme. It should be stressed that the MDE-based approach for implementing a derivative-free EKF technique demands the Cholesky factorization in each iterate of the auxiliary function for computing the sample points required by the right-hand side expression in~\eqref{eq2.4b}. This makes the MDE-based implementation strategy additionally vulnerable
to roundoff. Again, the filtering process is interrupted when the Cholesky factorization is unfeasible. We anticipate that the numerical stability of the MDE-based derivative-free EKF implementation is noticeably worse compared to the SPDE-based derivative-free EKF framework.

To conclude this section, we note that the MDEs in~\eqref{eq2.4a}, \eqref{eq2.4b} consist of $n^2+n$ equations as well as the SPDEs in~\eqref{SPDEs:eqs} contain $n^2$ equations to be solved together with the mean equation of size $n$. Thus, the MDE- and SPDE-based derivative-free EKF implementation ways outlined in Algorithms~1 and~2, respectively, are of the same computational cost, but the SPDE-based Algorithm~2 is expected to be more numerically stable than the MDE-based counterpart in Algorithm~1.

\section{Derivation of Cholesky factored-form continuous-discrete derivative-free EKF methods} \label{sec:main:squareroot}

The square-root (SR) filtering methods are known to improve a numerical robustness of filters' implementations with respect to roundoff errors. Additionally, they ensure the theoretical properties of the filter covariance matrix in a finite precision arithmetics, i.e. its symmetric form and positive definiteness. A variety of the SR filtering algorithms comes from the chosen factorization of the form $P=SS^{\top}$. In this paper, we design the Cholesky factored-form derivative-free EKF methods, i.e. we have the lower triangular SR matrix $S_{k|k}:=P_{k|k}^{1/2}$ at the measurement update step as well as the lower triangular SR factor  $S_{k|k-1}:=P_{k|k-1}^{1/2}$ at the prediction filtering step. As mentioned previously, the derivation of the related SR implementations is of special interest for all KF-like sampled-data estimators (e.g. CKF, GHQF, UKF and derivative-free EKF), because they demand the Cholesky factorization in each iterate of the filtering recursion that can be infeasible; see also the discussion in~\cite[p.~1262]{Haykin2009}. Unlike the {\it conventional} implementation way, the SR filtering methodology requires only one Cholesky decomposition, which is performed at the initial filtering step for the given initial matrix $\Pi_0 > 0$. Next, the filtering equations under examination should be re-derived in terms of propagating and updating the Cholesky factors of covariance matrices, only. This routine ensures the positive (semi-) definiteness and symmetric form of the resulted error covariance $P_{k|k}^{1/2}P_{k|k}^{\top/2} = P_{k|k}$ at any time instance $t_k$ and yields the improved numerical stability to roundoff errors; see more detail in~\cite[Chapter~7]{2015:Grewal:book}.

In summary, the goal of this section is to express the novel MDE- and SPDE-based MATLAB-oriented continuous-discrete derivative-free EKF techniques outlined in Table~\ref{Tab:1} in terms of the Cholesky factors of the filter covariance matrices. We start with the MDE-based derivative-free EKF approach. Taking into account the result proved in Lemma~\ref{Lm3-1}, we formulate the Cholesky-based SR MDEs to be solved at the prediction step of the continuous-discrete derivative-free EKF over the sampling interval $[t_{k-1}, t_k]$:
\begin{align}
\hat{x}'(t) & =  {f}\bigl(t, \hat x(t)\bigr), \label{sr:eq2.4a} \\
 \bigl[P^{1/2}(t)\bigr]' & = P^{1/2}(t) \; \Phi\Bigl( P^{-1/2}(t) \; M(t) \; P^{-\top/2}(t)\Bigr) \label{sr:eq2.4b} \\
\mbox{ with } \hat x(t_{k-1}) & = \hat x_{k-1|k-1},\quad P^{1/2}(t_{k-1})=P^{1/2}_{k-1|k-1}. \label{sr:eq2.4c}
\end{align}
The matrix $M(t)$ and operator $\Phi\bigl(A(t)\bigr)$ are defined in formulas~\eqref{eq3.9} and~\eqref{eq3.10}, i.e.
\begin{align*}
M(t) &:= \frac{\alpha}{\sqrt{n}} \left[ P^{1/2}(t) \; \overline{{\mathbb F}{\mathbb X}}^{\top}(t) + \overline{{\mathbb F}{\mathbb X}}(t) P^{\top/2}(t)\right] + G(t)\,Q(t)\,G^\top\!(t)\\
\Phi\bigl(A(t)\bigr) & :=\overline{L}(t) + 0.5\,{D}(t).
\end{align*}

Next, it is worth noting here that the prediction step of the SPDE-based derivative-free EKF summarized in Algorithm~2 works with the sample points, directly. It propagates them through equation~\eqref{SPDEs:eqs}. Hence, the SPDE-based prediction step does not require any changes to design its SR SPDE-based counterpart, i.e. it remains the same as in Algorithm~2.

Let us consider, the measurement update step of the novel MDE- and SPDE-based Matlab-oriented continuous-discrete derivative-free EKF techniques outlined in Table~\ref{Tab:1}. As can be seen, the filtering equations in~\eqref{ckf:state}, \eqref{ckf:gain} and~\eqref{eq:approx:Pxz:new} should be written in terms of the SR covariance matrix $P^{1/2}_{k|k}$. Recall, the term $P^{1/2}_{k|k-1}$ is calculated at the previous prediction step.

It is well known that the SR filters' implementations with the use of orthogonal transformations provide an additional numerical stability in a finite precision arithmetics. There, the orthogonal transformations are used for computing the Cholesky factorization of a positive definite matrix that obeys the equation $C = AA^{\top} + BB^{\top}$ by applying the orthogonal transformation to the pre-array $[A \; B]$ as follows: $[A \; B] Q = [R \; 0]$ where $R$ is a lower triangular Cholesky factor of the matrix $C$ as we are looking for.

Let us consider the first equation in~\eqref{eq:approx:Pxz:new}. We factorize it in the following form:
\begin{align}
R_{e,k} & = \overline{\mathbb Z}_{k|k-1}\overline{\mathbb Z}_{k|k-1}^{\top} + R_k = \left[\overline{\mathbb Z}_{k|k-1} \;\; R_k^{1/2}\right]\left[\overline{\mathbb Z}_{k|k-1} \;\; R_k^{1/2}\right]^{\top}
\label{eq:proof:4}
\end{align}
where $R_k^{1/2}$ is the lower triangular Cholesky factor of the measurement noise covariance $R_k$, i.e. $R_k = R_k^{1/2}R_k^{\top/2}$.

Thus, we have an equality
\[\left[\overline{\mathbb Z}_{k|k-1} \quad R_k^{1/2}\right]Q = [R_{e,k}^{1/2} \quad 0] \]
where $Q$ is any orthogonal matrix that lower triangulates the pre-array, i.e. $R_{e,k}^{1/2}$ is the lower triangular Cholesky factor that we are looking for.

Having substituted the term $R_{e,k}^{1/2}$ into the first equation in~\eqref{ckf:gain}, we obtain the formula for the filter's gain computation in terms of the SR error covariance factors as follows:
\[ K_{k} = P_{xz,k}R_{e,k}^{-1} = P_{xz,k}R_{e,k}^{-\top/2}R_{e,k}^{-1/2}. \]

Finally, we need to factorize the second equation in~\eqref{ckf:gain}. Unfortunately, it has the form of $C = AA^{\top} - BB^{\top}$ and the usual orthogonal transformations can not be applied for calculating the SR factor of $P_{k|k}$. However, we can derive a symmetric equation for calculating $P_{k|k}$ that can be further factorized in usual way. Taking into account that $K_k = P_{xz,k}R_{e,k}^{-1}$ and the symmetric form of the covariance matrix $R_{e,k}$, we prove
\begin{equation}
{K}_{k}R_{e,k}{K}_{k}^{\top} = P_{xz,k}{K}_{k}^{\top} \quad \mbox{ and } \quad {K}_{k}R_{e,k}{K}_{k}^{\top} = {K}_{k}P_{xz,k}^{\top} \label{ppppr:1}
\end{equation}
Taking into account $P_{xz,k}  = \overline{\mathbb  X}_{k|k-1} \overline{\mathbb  Z}_{k|k-1}^\top$, $P_{k|k-1} = \overline{\mathbb  X}_{k|k-1} \overline{\mathbb  X}_{k|k-1}^\top$,  $R_{e,k} = \overline{\mathbb Z}_{k|k-1}\overline{\mathbb Z}_{k|k-1}^{\top} + R_k$ and having substituted formula~\eqref{ppppr:1} into  the second equation in~\eqref{ckf:gain}, we derive
\begin{align}
P_{k|k} &= P_{k|k-1} - K_k R_{e,k}K_k^{\top}  = P_{k|k-1} - {K}_{k}P_{xz,k}^{\top} + K_k^{\top} R_{e,k}K_k  \nonumber \\
& - K_k^{\top} R_{e,k}K_k = P_{k|k-1} - P_{xz,k}{K}_{k}^{\top} + K_k R_{e,k}K_k^{\top} - {K}_{k}P_{xz,k}^{\top} \nonumber \\
& = \overline{\mathbb  X}_{k|k-1} \overline{\mathbb  X}_{k|k-1}^\top + \overline{\mathbb  X}_{k|k-1} \overline{\mathbb  Z}_{k|k-1}^\top K_k^\top + K_k (\overline{\mathbb Z}_{k|k-1}\overline{\mathbb Z}_{k|k-1}^{\top} + R_k) K_k^{\top} \nonumber \\
& - K_k  \overline{\mathbb  Z}_{k|k-1} \overline{\mathbb  X}_{k|k-1}^\top = K_k R_k K_k^{\top} \nonumber \\
& + \left[\overline{\mathbb X}_{k|k-1}-{K}_{k}\overline{\mathbb Z}_{k|k-1}\right]\left[\overline{\mathbb X}_{k|k-1}-{K}_{k}\overline{\mathbb Z}_{k|k-1}\right]^{\top}.
\end{align}

Thus, we get an equality
\[\left[ (\overline{\mathbb X}_{k|k-1}-{K}_{k}\overline{\mathbb Z}_{k|k-1}) \quad K_kR_k^{1/2}\right]Q = [P_{k|k}^{1/2} \quad 0] \]
where $Q$ is any orthogonal matrix that lower triangulates the pre-array, i.e. $P_{k|k}^{1/2}$ is the lower triangular Cholesky factor that we are looking for.

Finally, we may suggest one more SR implementation method by taking into account that equations~\eqref{ckf:gain} and~\eqref{eq:approx:Pxz:new} can be summarized in the following equality:
\[
\underbrace{\begin{bmatrix}
\overline{\mathbb  Z}_{k|k-1}  & R_k^{1/2} \\
\overline{\mathbb  X}_{k|k-1}  & 0
\end{bmatrix}}_{pre-array} Q =
\underbrace{\begin{bmatrix}
R_{e,k}^{1/2} & 0\\
\bar P_{xz,k} & P^{1/2}_{k|k}
\end{bmatrix}}_{post-array}
\]
where $Q$ is any orthogonal matrix that (block) lower triangulates the pre-array, i.e. $R_{e,k}^{1/2}$ and $P_{k|k}^{1/2}$ are the lower triangular Cholesky factors, which we are looking for. They are simply read-off from the post-array.

\begin{table*}[ht!]
{\small
\renewcommand{\arraystretch}{1.3}
\caption{The {\it square-root} continuous-discrete MATLAB-oriented derivative-free EKF implementation methods.} \label{Tab:2}
\centering
\begin{tabular}{l||l|l}
\hline
& \cellcolor{myGray} {\bf SR MDE-based derivative-free EKF: Algorithm~1a} & \cellcolor{myGray} {\bf SR SPDE-based derivative-free EKF: Algorithm~2a} \\
\hline
\hline
\textsc{Initialization:}  &  \multicolumn{2}{l}{0. Cholesky decomposition $\Pi_0=\Pi_0^{1/2}\Pi_0^{\top/2}$. Set $\hat x_{0|0} = \bar x_0$, $P_{0|0}^{1/2} = \Pi_0^{1/2}$ as well as $\alpha$ and the ODE solver's options in~\eqref{eq2.136}.  }\\
\hline
\textsc{Time} & 1. Form matrix $XP_{k-1|k-1} = [\hat x_{k-1|k-1}, P^{1/2}_{k-1|k-1}]$. & 1. Generate ${\mathbb  X}_{k-1|k-1}=\hat x_{k-1|k-1} {\mathbf 1}^{\top} + \frac{\sqrt{n}}{\alpha}P_{k-1|k-1}^{1/2}$. \\
\textsc{Update (TU):} & 2. Reshape $x^{(0)}_{k-1} = XP_{k-1|k-1}\verb"(:)"$. & \qquad $\rightarrow$ Repeat steps~2-6 from Algorithm~2. \\
& 3. Integrate $x_{k|k-1}\leftarrow \texttt{odesolver[SR MDEs},x^{(0)}_{k-1},[t_{k-1},t_k]]$. &  \\
& 4. $XP_{k|k-1} \leftarrow \texttt{reshape}(x_{k|k-1}^{\texttt{end}},n,n+1)$.  &  \\
& 5. Recover $\hat x_{k|k-1} = [XP_{k|k-1}]_1$. &  \\
& 6. Recover $P_{k|k-1}^{1/2} = [XP_{k|k-1}]_{2~:~n+1}$. &  \\
\hline
\textsc{Measurement}  & 7. Define all sample points ${\mathbb  X}_{k|k-1}=\hat x_{k|k-1} {\mathbf 1}^{\top} + \frac{\sqrt{n}}{\alpha}P_{k|k-1}^{1/2}$. & $-$ The matrix ${\mathbb  X}_{k|k-1}$ is already available from TU. \\
\cline{2-3}
\textsc{Update (MU):} & \multicolumn{2}{l}{8. Propagate all sample points ${\mathcal X}_{i,k|k-1}$ from ${\mathbb  X}_{k|k-1}$: ${Z}_{i,k|k-1} = h\left(k,{\mathcal X}_{i,k|k-1}\right)$, $i=1,\ldots,n$, and the estimate $\hat z_{k|k-1} = h\left(k,\hat x_{k|k-1}\right)$. } \\
& \multicolumn{2}{l}{9. Define $\overline{\mathbb  X}_{k|k-1}  = \frac{\alpha}{\sqrt{n}} \bigl[{\mathcal X}_{1,k|k-1} - \hat x_{k|k-1} | \ldots | {\mathcal X}_{n,k|k-1} - \hat x_{k|k-1}\bigr]$ and $\overline{\mathbb  Z}_{k|k-1}  = \frac{\alpha}{\sqrt{n}}\bigl[{Z}_{1,k|k-1} - \hat z_{k|k-1}| \ldots | {Z}_{n,k|k-1} - \hat z_{k|k-1}\bigr]$.  } \\
& \multicolumn{2}{l}{10. Build pre-array $A = \left[\overline{\mathbb Z}_{k|k-1} \quad R_k^{1/2}\right]$. Compute $R = qr(A^{\top})$ where $R^{\top} = [R_{e,k}^{1/2} \quad 0]$. Read-off $R_{e,k}^{1/2}$, which we are looking for. } \\
& \multicolumn{2}{l}{11. Compute $P_{xz,k}=\overline{\mathbb X}_{k|k-1}\overline{\mathbb Z}_{k|k-1}^{\top}$, ${K}_{k}=P_{xz,k}R_{e,k}^{-\top/2}R_{e,k}^{-1/2}$.  Find $\hat x_{k|k}=\hat x_{k|k-1}+{K}_k(z_k-\hat z_{k|k-1})$.} \\
& \multicolumn{2}{l}{12. Build pre-array $A = \left[ (\overline{\mathbb X}_{k|k-1}-{K}_{k}\overline{\mathbb Z}_{k|k-1}) \quad K_kR_k^{1/2}\right]$. Compute $R = qr(A^{\top})$ where $R^{\top} = [P_{k|k}^{1/2} \quad 0]$. Read-off $P_{k|k}^{1/2}$. } \\
\hline
\hline
& \cellcolor{myGray} {\bf SR MDE-based derivative-free EKF: Algorithm~1b} & \cellcolor{myGray} {\bf SR SPDE-based derivative-free EKF: Algorithm~2b} \\
\hline
\hline
\textsc{Initialization:}  &  \qquad $\rightarrow$ Repeat from Algorithm~1a. &  \qquad $\rightarrow$ Repeat from Algorithm~2a.  \\
\textsc{Time Update} & \qquad $\rightarrow$ Repeat from Algorithm~1a. &  \qquad $\rightarrow$ Repeat from Algorithm~2a.  \\
\hline
\textsc{Measurement}  & 7. Define all sample points ${\mathbb  X}_{k|k-1}=\hat x_{k|k-1} {\mathbf 1}^{\top} + \frac{\sqrt{n}}{\alpha}P_{k|k-1}^{1/2}$. & $-$ The matrix ${\mathbb  X}_{k|k-1}$ is already available from TU. \\
\cline{2-3}
\textsc{Update (MU):} & \multicolumn{2}{l}{8. Propagate all sample points ${\mathcal X}_{i,k|k-1}$ from ${\mathbb  X}_{k|k-1}$: ${Z}_{i,k|k-1} = h\left(k,{\mathcal X}_{i,k|k-1}\right)$, $i=1,\ldots,n$, and the estimate $\hat z_{k|k-1} = h\left(k,\hat x_{k|k-1}\right)$. } \\
& \multicolumn{2}{l}{9. Define $\overline{\mathbb  X}_{k|k-1}  = \frac{\alpha}{\sqrt{n}} \bigl[{\mathcal X}_{1,k|k-1} - \hat x_{k|k-1} | \ldots | {\mathcal X}_{n,k|k-1} - \hat x_{k|k-1}\bigr]$ and $\overline{\mathbb  Z}_{k|k-1}  = \frac{\alpha}{\sqrt{n}}\bigl[{Z}_{1,k|k-1} - \hat z_{k|k-1}| \ldots | {Z}_{n,k|k-1} - \hat z_{k|k-1}\bigr]$.  } \\
& \multicolumn{2}{l}{10. Build pre-array $A = \begin{bmatrix}
\overline{\mathbb  Z}_{k|k-1}  & R_k^{1/2} \\
\overline{\mathbb  X}_{k|k-1}  & 0
\end{bmatrix}$. Compute $R = qr(A^{\top})$ where $R^{\top} = \begin{bmatrix}
R_{e,k}^{1/2} & 0\\
\bar P_{xz,k} & P^{1/2}_{k|k}
\end{bmatrix}$. Read-off $R_{e,k}^{1/2}$, $\bar P_{xz,k}$, $P_{k|k}^{1/2}$. } \\
& \multicolumn{2}{l}{11. Compute $K_k = \bar P_{xz,k} R_{e,k}^{-1/2}$. Find $\hat x_{k|k}=\hat x_{k|k-1}+{K}_k(z_k-\hat z_{k|k-1})$.} \\
\hline
\hline
Auxiliary & $[\tilde x(t)] \leftarrow \proc{SR MDEs}(x(t),t,n,\alpha,G,Q)$
&
$[\widetilde{{\mathbb  X}}(t)] \leftarrow \proc{SPDEs}(\widetilde{{\mathbb  X}}(t),t,n,\alpha,G,Q)$\\
Functions & Get matrix $X = \verb"reshape"(x(t),n,n+1)$; & \qquad $\rightarrow$ Repeat from Algorithm~2. \\
& Recover $\hat x(t) = [X]_1$, $P^{1/2}(t) = [X]_{2~:~n+1}$. Define sample points~\eqref{SP:equation:1}; & \\
& Propagate $d{\hat x}(t)/dt:= f\bigl(t,{\hat x}(t)\bigr)$ and all sample points to find $\overline{\mathbb FX}(t)$; & \\
& Find $M(t)$ by~\eqref{eq3.9} and split $M=\bar L + D + \bar U$. Find $\Phi(M) = \bar L + 0.5 D$; & \\
& Compute the right-hand side of the SR MDEs in~\eqref{sr:eq2.4b}; & \\
& Collect extended matrix $\tilde X = [d\hat x(t)/dt, dP^{1/2}(t)/dt]$; & \\
& Reshape into a vector form $\tilde x(t)=\tilde X\verb"(:)"$. & \\
\hline
\end{tabular}
}
\end{table*}

The equality above can be proved by multiplying the pre- and post-arrays involved. Having been compared both sides of the obtained formulas with  equations~\eqref{ckf:gain} and~\eqref{eq:approx:Pxz:new}, the proof is straightforward. Besides, the term $\bar P_{xz,k} :=P_{xz,k}R_{e,k}^{-{\top}/2} = \overline{\mathbb  X}_{k|k-1} \overline{\mathbb  Z}_{k|k-1}^\top R_{e,k}^{-{\top}/2}$ and, hence, the gain matrix is calculated in terms of SR factors $R_{e,k}^{1/2}$ as follows: $K_k = P_{xz,k} R_{e,k}^{-1} = \bar P_{xz,k} R_{e,k}^{-1/2}$.

To summarize, we suggest two variants of the SR derivative-free EKF measurement update. Thus, we obtain the Cholesky-based MDE- and SPDE-based derivative-free EKF methods derived above in the form of Algorithms~1a, 1b and Algorithms~2a, 2b, respectively. They are presented in Table~\ref{Tab:2}. It should be stressed that the QR factorization is performed by the built-in MATLAB function \texttt{qr(A)} that yields an upper triangular matrix, i.e. $Q$ is an orthogonal rotation that upper triangulates the given pre-arrays. In the MATLAB-oriented SR algorithms in Table~\ref{Tab:2} this fact is taken into account in order to formulate the methods in terms of lower triangular Cholesky factors.

Let us discuss the novel SR MDE- and SPDE-based MATLAB-oriented continuous-discrete derivative-free EKF techniques summarized in Table~\ref{Tab:2}.
All algorithms are of {\it covariance SR}-type implementations because the Cholesky factors of the filter error covariance matrix are propagated and updated. The Cholesky decomposition is avoided at all filtering steps, except the initial stage. Besides, all SR algorithms proposed utilize numerically stable orthogonal transformations as far as possible for propagating the Cholesky factors. This additionally improves the numerical stability of the suggested SR implementations. We also note that SR Algorithms~1a and~2a require two QR factorizations at the measurement update step. Meanwhile SR Algorithms~1b and~2b demand only one QR transformation at each iterate, i.e. they are a bit faster than their counterparts in Algorithms~1b and~2b. Finally, we compare the SR MDE-based implementation in Algorithm~1a and the related SR SPDE-based solution in Algorithm~2a as well as Algorithm~1b with Algorithm~2b. We conclude that the SR MDE- and the related SR SPDE-based solutions consist of the same number of equations to be solved at the prediction step and include the same measurement updates, i.e. we conclude that Algorithms~1a and~2a are of the same cost as well as Algorithms~1b and~2b possess the same computational expenses.

In summary, the SR derivative-free EKF implementation methods derived in this section are more numerically stable with respect to roundoff errors but more computationally expensive due to QR factorizations compared to the conventional Algorithms~1 and~2 in Table~\ref{Tab:1}.

\section{Numerical experiments} \label{numerical:experiments}

To illustrate the performance of the novel continuous-discrete derivative-free EKF methods derived, we provide a set of numerical tests. Following the suggestion in~\cite{quine2006derivative}, we implement all derivative-free EKF methods with $\alpha = 10^{3}$. Recall, the novel MATLAB-oriented continuous-discrete EKF methods provide a flexibility  in a choice of the built-in ODE solver to be used. In the first numerical tests, we utilize the built-in ODE solver \verb"ode45" with the tolerance value $\epsilon_g = 10^{-4}$ in all novel MDE- and SPDE-based EKF algorithms. More precisely, the \texttt{ode45} employs the explicit embedded Runge-Kutta pair of orders~4 and~5 with the scaled local error control designed in~\cite{DoPr80}. It is accepted as a benchmark code for integrating nonstiff ODE models in practice. Additionally, we utilize the standard continuous-discrete EKF with the same ODE solver and the tolerance value of $\epsilon_g = 10^{-4}$ in order to provide a fair comparative study with the novel methods; see the summary of the implementation method in~\cite[Algorithm~A.2]{KuKu16SISCI}.

\begin{figure*}
\begin{tabular}{cc}
\includegraphics[width=0.5\textwidth]{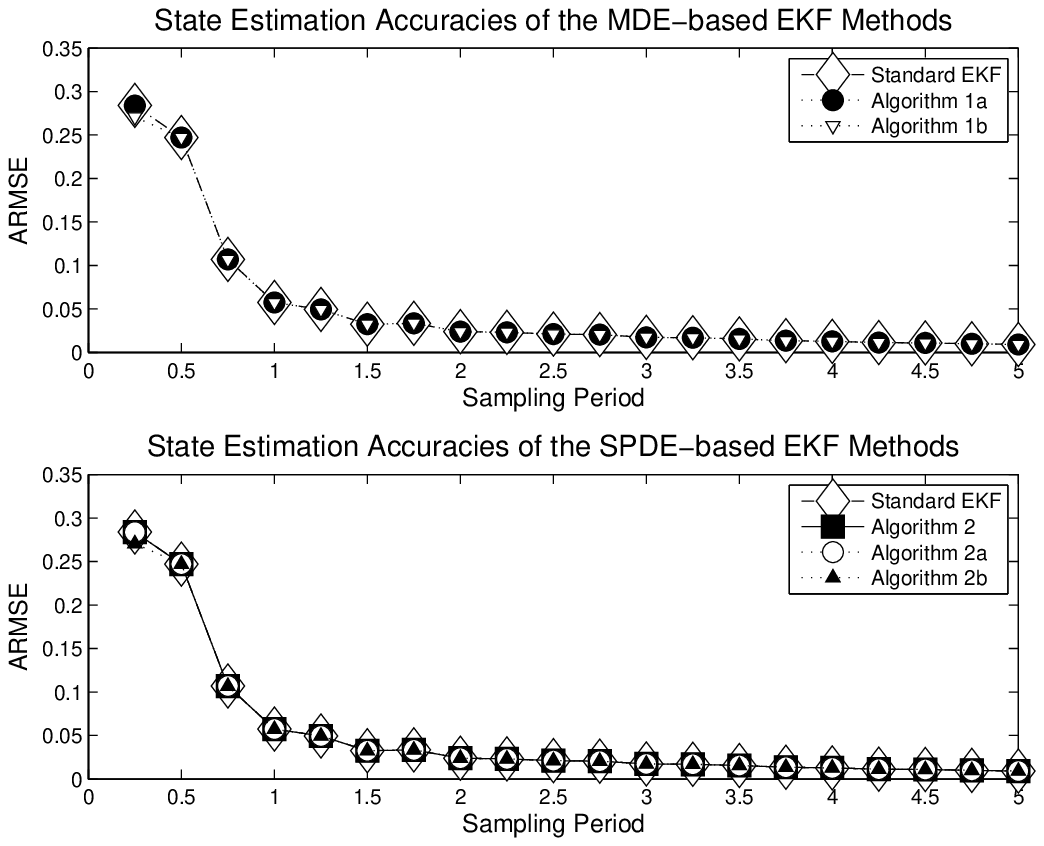} &  \includegraphics[width=0.5\textwidth]{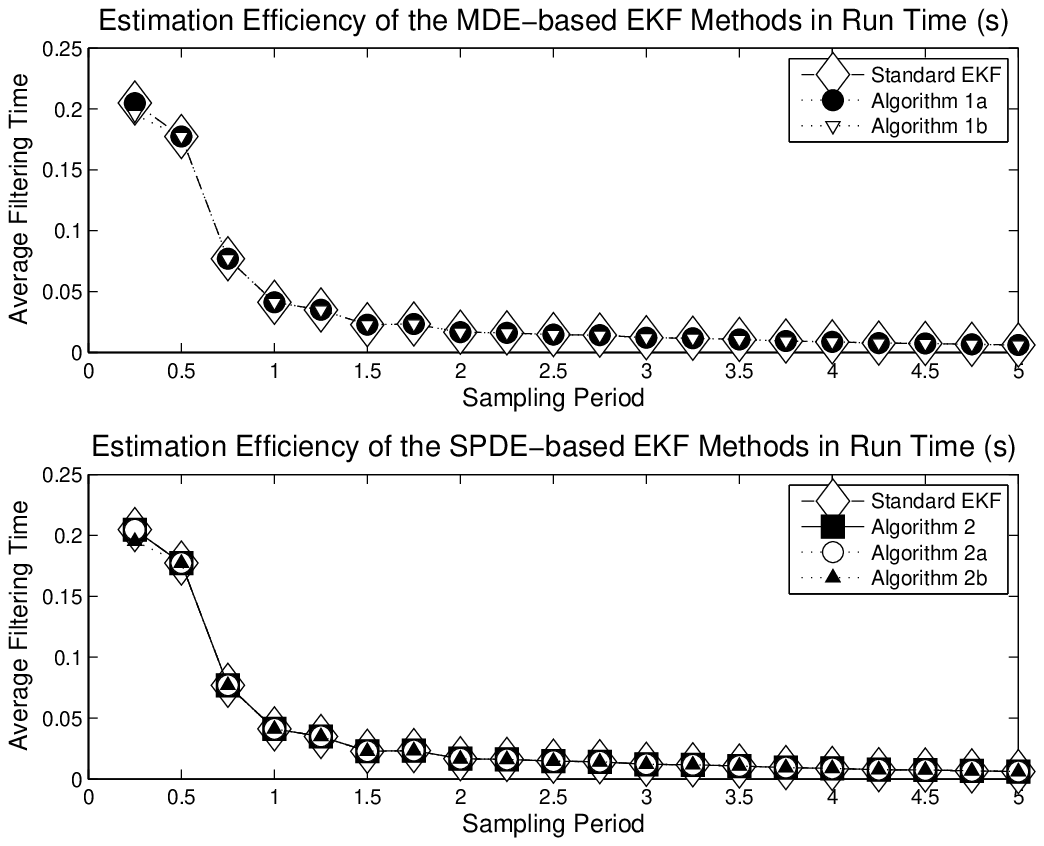}
\end{tabular}
\caption{The accuracies (left graph) and efficiencies (right graph) of various continuous-discrete EKF estimators on the test problem in Example~\ref{ex:2}.} \label{fig:1:new}
\end{figure*}

\begin{exmp}[Gas-phase reversible reaction in CSTR] \label{ex:2}
Following~\cite{KuKu19IJRNC}, the gas-phase reversible reaction with three species denoted as $A$, $B$ and $C$ are given as follows:
\begin{equation}\label{eq2.1}
A\quad {{k_1\atop\rightleftharpoons}\atop {\scriptstyle k_2}}\quad B+C,\quad
2B\quad {{k_3\atop\rightleftharpoons}\atop {\scriptstyle k_4}}\quad B+C,
\end{equation}
where the fixed coefficients $k_1=0.5$, $k_2=0.05$, $k_3=0.2$ and $k_4=0.01$. The stoichiometric matrix $\nu$ of reaction (\ref{eq2.1}) and the reaction rates $r$ are chosen to be
\begin{equation}\label{eq2.2}
\nu=\left[
\begin{array}{ccc}
-1 & 1 & 1\\
0 & -2 & 1
\end{array}
\right], \quad
r=\left[
\begin{array}{c}
k_1 c_A  -k_2 c_B c_C\\
k_3 c_B^2 -k_4 c_C
\end{array}
\right],
\end{equation}
in which $c_A$ stands for the concentration of the species $A$ in moles per liter, and so on. Thus, the state of this chemical system is defined by the vector $x(t)=\left[c_A\quad c_B\quad c_C\right]^\top$.

We explore reaction~\eqref{eq2.1}, \eqref{eq2.2} in a continuously stirred tank reactor (CSTR). The well-mixed, isothermal CSTR is simulated by the following SDE:
\begin{equation}\label{eq2.5:new}
dx(t)=\left[\frac{Q_f}{V_R} c_f - \frac{Q_0}{V_R} x(t) + \nu^\top r\right] dt+G dw(t), \quad t>0,
\end{equation}
where $c_f=x_0$, the coefficients $Q_f=Q_0=1$, $V_R=100$, and the matrix $\nu$ and the vector $r$ are given in (\ref{eq2.2}). The continuous-time process noise $dw(t)$ is the zero-mean white Gaussian process with the diagonal covariance matrix $Q\,dt=\mbox{\rm diag}\{10^{-6}/\delta,10^{-6}/\delta,10^{-6}/\delta\}\,dt$, and the constant diffusion matrix is $G=I_3$. The initials are $\bar x_0 = [0.5, 0.05, 0]^{\top}$ and $\Pi_0 = I_3$.

The measurement equation is given in the following form:
\begin{equation}\label{eq4.5:new}
z_k=\left[RT\quad RT\quad RT\right] x_k + v_k, \quad
  \begin{array}{l}
    v_k \sim {\cal N}(0,R_z); \\
    R_z  = 0.25^2
\end{array}
\end{equation}
where $R$ is the ideal gas constant and $T$ is the reactor temperature in Kelvin, i.e. $RT=32.84$.
\end{exmp}

To assess the estimation accuracies of the novel estimators developed in this paper, we first solve the direct problem that is the numerical simulation of the given model. More precisely, we solve the SDE in~\eqref{eq2.5:new} with a small stepsize $\delta_t = 10^{-3}$ on interval $[0, 30]$(s) to generate the {\it true} state vector $x^{\rm true}(t_k)$, $t_k \in [0, 30]$(s). Next, the measurement model~\eqref{eq4.5:new} is utilized for creating the history of simulated measurements with various sampling rates $\Delta =0.5, 1,\ldots, 4.5, 5$(s). Finally, for each fixed $\Delta$(s) value, the inverse problem, i.e. the filtering problem, is solved to get the estimated hidden state $\hat x_{k|k}$ over the time interval interval $[0, 30]$(s). We repeat the numerical test for $100$ trials and compute the accumulated root mean square error (ARMSE) by averaging over $100$ Monte Carlo runs and three entries of the state vector $x(t) = \left[c_A\quad c_B\quad c_C\right]^\top$ as follows:
 \begin{align}
\mbox{\rm ARMSE} & =\Bigl[\frac{1}{M
K}\sum_{M=1}^{100}\sum_{k=1}^K\sum_{j=1}^{n}\bigl(x^{\rm true}_{k,j}-\hat
x_{k|k,j}\bigr)^2\Bigr]^{1/2} \label{eq:acc:1}
\end{align}
where the subindex $j$, $j =1,\ldots,n$, refers to the $j$th entry of the $n$-dimensional state vector.

In each Monte Carlo run, we obtain its own simulated ``true'' state trajectory and the measurement data history. It is important that all estimators utilize the same initial conditions, the same simulated ``true'' state trajectory and the same measurements. The ARMSE values computed for each sampling period $\Delta$ mentioned above by using the MDEs-based EKF methods (top graph) and the SPDE-based EKF implementations (bottom graph) are illustrated by the left-hand side of Fig.~\ref{fig:1:new}.  We also compute the average CPU time in seconds over $100$ Monte Carlo runs for each estimator under discussion and illustrate the resulted values by the right-hand side graphs of Fig.~\ref{fig:1:new}.

Having analyzed the results illustrate by  Fig.~\ref{fig:1:new}, we make a few conclusions. Firstly, we note that the conventional MDE-based Algorithm~1 fails to solve the filtering problem due to unfeasible Cholesky decomposition. This result was anticipated and discussed previously. This implementation method  demands the Cholesky factorization in each iterate of the auxiliary function for computing the sample points required by the right-hand side expression in~\eqref{eq2.4b}. This makes the conventional MDE-based approach highly  vulnerable to roundoff. Meanwhile, the SR MDE-based counterparts in Algorithms~1a and~1b easily solve the stated filtering problem. Hence, the SR approach indeed improves the numerical stability of the MDE-based continuous-discrete EKF technique.

We also observe that the SPDE-based frameworks, both the conventional and the square-root one, accurately solve the filtering problem at hand. As it was anticipated, the conventional SPDE-based Algorithm~2 is more stable to roundoff than the conventional MDE-based Algorithm~1, meanwhile, the SR MDE- and SR SPDE-based Algorithms~1a, 1b and~2a, 2b are all of improved numerical stability.

Secondly, all novel derivative-free EKF implementation methods, except Algorithm~1, work accurately for any sampling interval length under examination due to the MATLAB built-in discretization error control involved. All of them treat the long sampling interval with the same estimation accuracies and similar CPU time. It is an useful feature of the novel estimators that is of special interest in case of missing measurements, i.e. in case of long and irregular sampling intervals.

Thirdly, we conclude that all MDE- and SPDE-based derivative-free EKF implementations, except Algorithm~1, provide the same estimation quality with the same computational costs. This fact additionally substantiates the mathematical derivation of the square-root  MDE- and SPDE-based derivative-free EKF methods. Besides, we observe a similar execution time for all novel derivative-free EKF methods; see the right-hand side graphs of Fig.~\ref{fig:1:new}.

Finally, let us compare the performances of the novel estimators developed in this paper and the standard continuous-discrete EKF method  in~\cite[Algorithm~A.2]{KuKu16SISCI}. Following Fig.~\ref{fig:1:new}, we conclude that all methods under examination behave in the same way. In other words, they maintain the same estimation accuracies at the same computational expenses. This result was anticipated because of the following reasons. Firstly, the simulation results provided in~\cite{quine2006derivative} have shown an excellent convergence of the derivative-free EKF to the standard EKF when $\alpha = 1000$. Thus, the results obtained in our numerical tests are in line with the conclusion made in~\cite{quine2006derivative}, i.e. the same conclusion is obtained for the {\it continuous-discrete} EKF variants when $\alpha = 1000$. Indeed, as can be seen from the left-hand side graphs of Fig.~\ref{fig:1:new} the novel derivative-free EKF methods and the standard continuous-discrete EKF technique ensure the same estimation quality. Secondly, the number of equations to be solved at the prediction step in the standard continuous-discrete EKF method is $n^2+n$, i.e. it coincides with the size of the MDE and SPDE systems to be integrated in the novel derivative-free EKF methods. Besides, the measurement update steps of the derivative-free EKF methods and the standard EKF are of similar computational costs because the filter covariance matrix of size $n$ is approximated by using $n$ vectors. The extra computational expenses of the derivative-free EKF methods compared to the standard EKF come from the Cholesky decompositions required for the sample vectors generation. Consequently, all methods under examination possess similar CPU times to run the codes.

\begin{figure}
\includegraphics[width=0.5\textwidth]{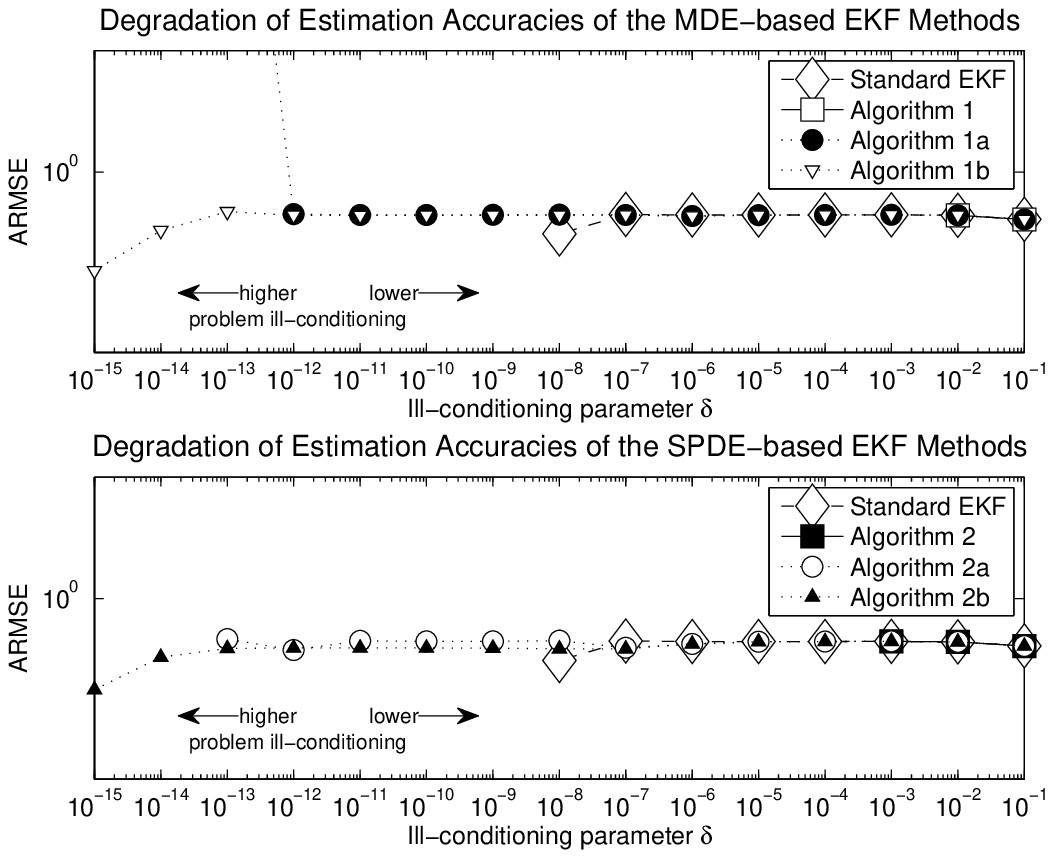}
\caption{The accuracies  degradation of various continuous-discrete EKF estimators on the ill-conditioned problem in Example~\ref{ex:2ill}.} \label{fig:2:new}
\end{figure}

Our next goal is to investigate a difference in the filters' numerical robustness with respect to roundoff errors. It is well known that the square-root implementation methods are more numerically stable than the conventional ones; see also the discussion in the previous section. To get more insights about the numerical robustness (to roundoff) of the novel  MDE- and SPDE-based derivative-free EKF algorithms, we follow the ill-conditioned test problem proposed in~\cite[Example~7.2]{2015:Grewal:book}. More precisely, a specially designed measurement scheme is suggested to be utilized for provoking the filters' numerical instability due to roundoff and for observing their divergence when the problem ill-conditioning increases.

\begin{exmp}[The CSTR ill-conditioned measurement scheme] \label{ex:2ill}
We design the following measurement scheme for provoking the filters' numerical instability due to roundoff:
\begin{align}
z_k & =
RT \begin{bmatrix}
1 & 1 & 1 \\
1 & 1 & 1 +\delta
\end{bmatrix}
x_k +
\begin{bmatrix}
v_k^1 \\
v_k^2
\end{bmatrix}, \;
  \begin{array}{l}
    v_k \sim {\cal N}(0,R_z); \\
     R_z=\delta^{2}I_2
\end{array} \label{measurementsill}
\end{align}
where parameter $\delta$ is used for simulating roundoff effect. This increasingly ill-conditioned target tracking scenario assumes that $\delta\to 0$, i.e. $\delta=10^{-1},10^{-2},\ldots,10^{-15}$.
\end{exmp}

For each value of the ill-conditioning parameter $\delta=10^{-1},10^{-2},\ldots,10^{-15}$, we repeat a set of numerical experiments explained above with the measurement model~\eqref{measurementsill} instead of model~\eqref{eq4.5:new}. Following~\cite[Example~7.2]{2015:Grewal:book}, the matrix $R_{e,k}$ becomes ill-conditioned while $\delta$ tends to machine precision limit. This destroys any KF-like implementation. We plot the ARMSE calculated for each filtering method under examination against the ill-conditioning parameter $\delta$ on Fig.~\ref{fig:2:new}.

Having analyzed the breakdown points illustrated by Fig.~\ref{fig:2:new}, we conclude that all conventional methods (which process the full filter covariance matrices), i.e. the MDE-based Algorithm~1, the SPDE-based Algorithm~2 and the standard EKF, degrade much faster than the square-root implementations while the ill-conditioning parameter $\delta$ tends to machine precision limit. This substantiates the improved numerical stability of the square-root methods derived. More precisely, the conventional MDE-based Algorithm~1 does not fail on well-conditioned tests in Example~\ref{ex:2ill} but it is again less stable than the SPDE-based Algorithm~2. Indeed, the conventional MDE-based Algorithms~1 fails for ill-conditioned state estimation problems with $\delta < 10^{-3}$ meanwhile the conventional SPDE-based Algorithm~2 solves the ill-conditioned tests at $\delta = 10^{-3}$ and fails for  $\delta < 10^{-4}$. Hence, we conclude that the MDE-based Algorithm~1 is more sensitive to roundoff errors than the SPDE-based Algorithm~2 as it was anticipated while the discussion in Section~\ref{sec:main:conventional}. It is interesting to note that the standard continuous-discrete EKF method in~\cite[Algorithm~A.2]{KuKu16SISCI} is also the conventional-type method but it seems  to be more robust than the derivative-free EKF counterparts since it fails for $\delta < 10^{-8}$. It seems that the accumulation of roundoff errors with the requirement to perform Cholesky decomposition in each iterate of the derivative-free Algorithms~1 and~2 worsen their numerical stability compared to the standard EKF framework.

Meanwhile, all square-root implementation methods easily manage the ill-conditioned tests at  $\delta = 10^{-4}$ and work accurately till $\delta = 10^{-13}$. Thus, we conclude that the square-root derivative-free EKF methods derived outperform their conventional counterparts for a numerical stability with respect to roundoff errors. Next, we observe that the SR MDE- and the SR SPDE-based derivative-free EKF Algorithms~1a and~2a fail when $\delta < 10^{-13}$ meanwhile the SR MDE- and the SR SPDE-based Algorithms~1b and~2b are able to solve these ill-conditioned estimation problems in accurate and robust way. Recall, Algorithms~1a and~2a require two QR factorizations on the measurement update step compared to Algorithms~1b and~2b, which demand only one QR transformation on each iterate. It seems that the accumulation of roundoff errors yields a worse performance of square-root Algorithms~1a and~2a compared to Algorithms~1b and~2b; see also the discussion in~\cite[Section 7.2.2]{2015:Grewal:book}. To summarize, we conclude that the square-root derivative-free EKF methods are more numerically stable than their conventional counterparts. Besides, Algorithms~1b and~2b are the most stable implementation methods among all square-root algorithms derived.

Our last numerical example is focused on exploring capacities of the novel filters for estimating stiff stochastic models. For that, we examine the stochastic Van der Pol oscillator, which is considered to be a classical benchmark in nonlinear filtering theory
by many authors~\cite{Frog2012,Ma08}. This test example can expose both nonstiff and stiff behaviors, depending on the value of its stiffness parameter $\lambda$. The difficulties of state estimation in stiff stochastic systems are observed by comparing performances of the filtering methods while $\lambda$ increases the period of oscillations. Our test example is rescaled as explained in~\cite[p.~5]{HaWa96}. Recall, one of the advantages of the suggested MATLAB-oriented derivative-free EKF implementation methods is their flexibility when any built-in ODE solver can be easily used for their implementation. Taking into account the potential stiffness of Example~\ref{ex:3}, we replace the \texttt{ode45} with \texttt{ode15s}, which is accepted as a benchmark code for integrating stiff ODE models in applied science and engineering. To perform a fair comparative study, the standard continuous-discrete EKF method in~\cite[Algorithm~A.2]{KuKu16SISCI} is also implemented by utilizing the MATLAB built-in ODE solver  \texttt{ode15s}.

\begin{figure}
\includegraphics[width=0.5\textwidth]{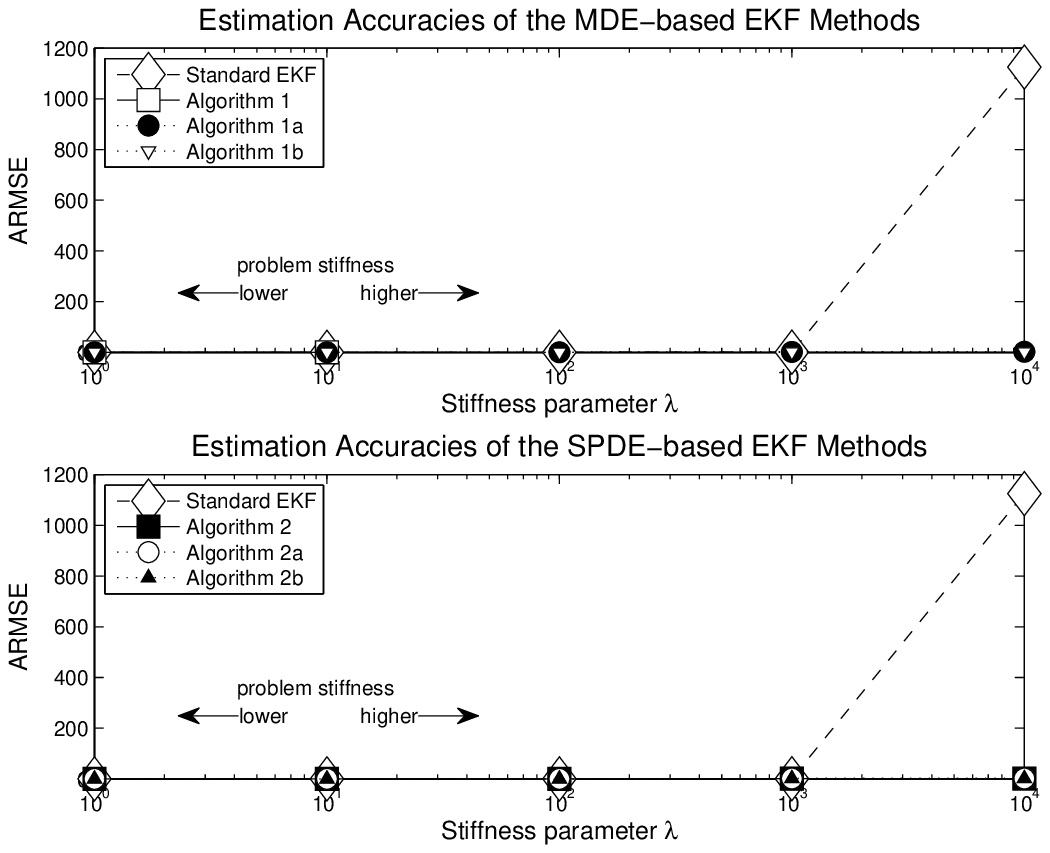}
\caption{The accuracies of various continuous-discrete EKF estimators depending on the value of the stiffness parameter $\lambda$ in Example~\ref{ex:3}.} \label{fig:3:new}
\end{figure}

\begin{exmp}[stochastic Van der Pol oscillator] \label{ex:3}
Consider the following SDE:
\begin{equation}\label{eq4.1:new2}
d\left[\!\!
\begin{array}{c}
x_1(t) \\
x_2(t)
\end{array}\!\!
\right]= \left[\!\!
\begin{array}{c}
x_2(t) \\
\lambda\bigl[(1-x_1^2(t))x_2(t)-x_1(t)\bigr]
\end{array}\!\!
\right]dt +\left[\!\!
\begin{array}{cc}
0 & 0 \\
0 & 1
\end{array}\!\!
\right]dw(t)
\end{equation}
where the initial state $\bar x_0=[2,0]^\top$ and $\Pi_0={\rm diag} \{0.1, 0.1 \}$ with process covariance $Q=I_2$.

The measurement equation is taken to be
\begin{equation}\label{eq4.2:new2} z_k =[1, \: 1] x_k + v_k, \quad   \begin{array}{l}
    v_k \sim {\cal N}(0,R); \\
    R  = 0.04.
\end{array}
\end{equation}
\end{exmp}
For each value of the stiffness parameter $\lambda = 10^{0}, 10^{1}, \ldots, 10^{4}$, we perform the following set of numerical experiments. The SDE in~\eqref{eq4.1:new2} is simulated with a small stepsize $\delta_t = 10^{-5}$ on interval $[0, 2]$(s) to generate the true state vector. Next, equation~\eqref{eq4.2:new2} is utilized for creating the history of simulated measurements with the sampling rate $\Delta = 0.2$(s). Given the measurement data, the filtering problem is solved and the  hidden state is estimated by the novel continuous-discrete derivative-free EKF methods. Again, we compute the ARMSE by averaging over $100$ Monte Carlo runs and plot them against the stiffness parameter $\lambda$ on Fig.~\ref{fig:3:new}.

Having analyzed the obtained results, we conclude that all continuous-discrete derivative-free EKF estimators suggested in this paper solve the nonstiff estimation problems with the same accuracies. Besides, the estimation quality is high because the estimation errors are small. The standard continuous-discrete EKF method provides the same estimation accuracies as the derivative-free EKF algorithms for the nonstiff estimation scenarios under examination.

However, the conventional MDE-based Algorithm~1 does not manage the stiff scenarios due to the unfeasible Cholesky decomposition required at any node of the adaptive mesh generated by the numerical integrator \texttt{ode15s} in use. As can be seen, the conventional MDE-based Algorithm~1 rapidly fails while the problem stiffness increases. More precisely, it fails for $\lambda > 10$. It is also important to note that the standard continuous-discrete EKF fails to solve the stiff problems with $\lambda > 10^{3}$ since the ARMSE is large. Meanwhile all other MATLAB-oriented derivative-free EKF implementation methods derived in this paper are able to solve the stiff estimation problems till $\lambda = 10^{4}$. In particular, this creates a solid background for using the novel square-root (both the MDE- and SPDE-based) continuous-discrete derivative-free EKF estimators for solving practical problems, including stiff scenarios as well as ill-conditioned cases.

\section{Concluding remarks}\label{Section:conclusion}

In this paper, we developed the {\it continuous-discrete} derivative-free EKF approach. For that, we derived the moment differential equations (MDEs) and sample point differential equations (SPDEs) for the derivative-free EKF prediction step. Thus, we obtained the MDE- and SPDE-based derivative-free EKF implementation frameworks, respectively. Additionally, the numerically stable square-root implementation methods are developed within the Cholesky decomposition of the filter covariance  matrices. More precisely, we derived the square-root version of the MDEs for the Cholesky factors propagation as well as several square-root implementations for the measurement update step. All novel square-root filtering algorithms utilize stable orthogonal transformations for updating the Cholesky factors of the filter covariance matrices as far as possible. Finally, the continuous-discrete derivative-free EKF methods developed are MATLAB-oriented implementations, i.e. any  built-in numerical integrator for solving either the MDEs or SPDEs can be used. The key benefit is that the numerical integration accuracy is controlled at prediction step in automatic way by using the local discretization error control that bounds the discretization errors occurred and makes the implementation methods accurate. This is done in automatic way and no extra coding is required from users. The results of numerical tests substantiate a high estimation quality and numerical robustness of the novel square-root (both the MDE- and SPDE-based) continuous-discrete derivative-free EKF estimators, including stiff problems.

Finally, it is worth noting here that the square-root implementation methods might be derived in alternative way that is by using the singular value decomposition (SVD) instead of the Cholesky factorization. The SVD-based filtering is still an open area for a future research. More precisely, the SVD solution for  the proposed  continuous-discrete derivative-free EKF framework requires the derivation of the MDEs and SPDEs in terms of the SVD factors of the filter error covariance matrix. This is an open problem for a future research.

\section*{Acknowledgements}
The authors acknowledge the financial support of the Portuguese FCT~--- \emph{Funda\c{c}\~ao para a Ci\^encia e a Tecnologia},
through the \emph{Scientific Employment Stimulus - 4th Edition} (CEEC-IND-4th edition) programme (grant number 2021.01450.CEECIND) and through the projects UIDB/04621/2020 and UIDP/04621/2020 of CEMAT/IST-ID, Center for Computational and Stochastic Mathematics, Instituto Superior T\'ecnico, University of Lisbon.

\section*{References}
\bibliographystyle{model1b-num-names}

\end{document}